\documentclass[10pt]{article}
\usepackage[T1]{fontenc} 
\usepackage[utf8]{inputenc}

\usepackage{lmodern}

\usepackage{secdot}
\usepackage{amssymb}
\usepackage{amsmath}
\usepackage{latexsym}
\usepackage{calrsfs}
\usepackage{graphicx,color,psfrag}
\usepackage{natbib}

\newtheorem{theorem}{Theorem}[section]
\newtheorem{lemma}[theorem]{Lemma}

\newtheorem{proposition}[theorem]{Proposition}

\newtheorem{assumption}[theorem]{Assumption}
\newenvironment{proof}{\trivlist\item[]\emph{Proof}.}
{\unskip\nobreak\hskip 1em plus 1fil\nobreak$\Box$
\parfillskip=0pt\endtrivlist}

\def\argmax{\mathop{\rm arg\,max}}

\newcommand{\Ex}{\ensuremath{\mathsf{E}}}

\usepackage[colorlinks=true,urlcolor=blue,citecolor=blue,linkcolor=blue]{hyperref}
\usepackage[noabbrev,nameinlink]{cleveref}     

\def\EMAIL#1{\href{mailto:#1}{#1}}

\begin{document}

\title{Computing a classic index for finite-horizon bandits}


\author{Jos\'e Ni\~no-Mora 
\\ Department of Statistics \\
    Carlos III University of Madrid \\
     28903 Getafe (Madrid), Spain \\  \EMAIL{jnimora@alum.mit.edu}, \href{http://alum.mit.edu/www/jnimora}{http://alum.mit.edu/www/jnimora} \\
      ORCID: \href{http://orcid.org/0000-0002-2172-3983}{0000-0002-2172-3983}}
 
\date{Published in \textit{INFORMS Journal on Computing}, vol.\  23, pp. 254--267,  2011 \\ \vspace{.1in}
DOI: \href{https://doi.org/10.1287/ijoc.1100.0398}{10.1287/ijoc.1100.0398}}

\maketitle


\begin{abstract}%
This paper considers the efficient exact computation of the counterpart of the Gittins index for a finite-horizon discrete-state bandit, which measures for each initial state the average productivity, given by the maximum ratio of expected total discounted reward earned to expected total discounted time expended that can be achieved through a number of successive plays stopping by the given horizon. Besides characterizing optimal policies for the finite-horizon one-armed bandit problem, such an index provides a suboptimal heuristic index rule for the intractable finite-horizon multiarmed bandit problem, which represents the natural extension of the Gittins index rule (optimal in the infinite-horizon case). Although such a finite-horizon index was introduced in classic work in the 1950s, investigation of its efficient exact computation has received scant attention. This paper introduces a recursive adaptive-greedy algorithm using only arithmetic operations that computes the index in (pseudo-)polynomial time in the problem parameters (number of project states and time horizon length). In the special case of a project with limited transitions per state, the complexity is either reduced or depends only on the length of the time horizon. The proposed algorithm is benchmarked in a computational study against the conventional calibration method.
\end{abstract}%

 
\newpage

\section{Introduction}
\label{s:intro}
This paper deals with a class of finite-horizon discrete-state bandit 
problems, whose optimal policy is known to be of index type. 
In contrast to the existing literature, where such an index is
computed
approximately via the so-called calibration method, this paper
provides an efficient and exact algorithm to compute the index.
\subsection{Finite-horizon multiarmed bandits}
\label{s:fhmab}
In the classic \emph{finite-horizon multiarmed bandit problem}
(FHMABP), a decision maker aims to maximize the expected total discounted
reward earned from
a finite collection of $M$ dynamic and stochastic projects,  one
 of which must be engaged at each of a finite number $T$ of discrete time 
periods  $t = 0, 1,
 \ldots, T-1$.
Project $m = 1, \ldots, M$ is modeled as a discrete-time \emph{bandit}, i.e., a
binary-action (active: $a_m(t) = 1$; passive: $a_m(t) = 0$)
\emph{Markov decision process} (MDP) whose state 
$X_m(t)$ moves through the discrete (finite or countably infinite)
state space $\mathbb{X}_m$.
If the project is 
engaged ($a_m(t) = 1$) at time  $t < T$  when it
occupies state $X_m(t) = i_m$,  it yields an expected active reward
$R_m^1(i_m) \equiv R_m(i_m)$  and its state evolves 
to $j_m$ with probability $p_m(i_m, j_m)$.
Otherwise, it neither yields reward (i.e., the passive reward
is $R_m^0(i_m) \equiv 0$) nor changes state.
Rewards are discounted with factor $0 < \beta \leq 1$, where
the term ``discounted'' is abused to include 
the  undiscounted case $\beta = 1$.

Decisions as to which  project  to engage at each time
 are based on
adoption of a \emph{scheduling policy} $\boldsymbol{\pi}$, to be drawn from
the
class $\boldsymbol{\Pi}$ of \emph{admissible policies}, which  engage one project at each time before time $T$, and are
nonanticipative (with respect to the history of elapsed time periods, states and actions) and possibly randomized.

The FHMABP
is to find an admissible policy that
maximizes the expected total discounted reward earned.
Denoting by 
$\Ex_{\mathbf{i}}^{\boldsymbol{\pi}}\left[\cdot\right]$ the
expectation  under policy $\boldsymbol{\pi}$ conditioned on the initial joint
state being equal to $\mathbf{i} = (i_m)$,
we can formulate such a problem as 
\begin{equation}
\label{eq:dsp}
\max_{\boldsymbol{\pi} \in \boldsymbol{\Pi}} 
\Ex_{\mathbf{i}}^{\boldsymbol{\pi}}\left[
\sum_{t=0}^{T-1} \sum_{m=1}^M   \beta^t R_m^{a_m(t)}\big(X_m(t)\big) \right].
\end{equation}

The problem has its roots in the seminal works of
 \cite{robbins52} and \cite{bjk56}, who focused
on the much-studied case where engaging a project corresponds to sampling from a 
Bernoulli population with unknown success probability, the goal being
to maximize the expected number of successes over $T$ plays.
An MDP formulation is obtained by a Bayesian
approach, where a project/population state is its posterior
distribution. 

The above FHMABP and some of its variants have since 
drawn extensive research attention, due  to their theoretical and
practical interest. See, e.g., the monograph by \cite{berryFrist85} and
references therein. 
More recently, \cite{carogall07} address a problem extension where 
 $K < M$ projects are to be engaged at each time, motivated
by a dynamic assortment problem in the fashion retail industry.

\subsection{The average-productivity index policy}
\label{s:acip}
Finding an optimal policy for such a problem through 
numerical solution of its \emph{dynamic programming} (DP) equations quickly
becomes computationally intractable as the time horizon or the
projects' state spaces grow, which has led researchers to
 investigate a variety of tractable, though
suboptimal,
heuristic scheduling rules. Simple examples include the ``\emph{play the
winner/switch from a loser}'' rule (for Bernoulli bandits) and the
\emph{myopic} policy, which engages at each time a project of currently highest
 expected reward.

Yet, for a special case of  (\ref{eq:dsp}), the 
\emph{two-armed bandit problem with one arm known} --- also known as the \emph{one-armed bandit problem}--- where
there are two projects and one (the \emph{known arm} or \emph{standard
project}) 
has a single state with
 reward $\lambda$, the structure of optimal policies
is well known, being characterized by an \emph{index}
$\lambda^*(d, i)$ attached to states $i \in \mathbb{X}$ and \emph{times-to-go}
$d = 1, \ldots, T$ for the other project (the \emph{unknown arm}, for which 
the label $m$ is henceforth dropped from the notation). 
Note that in such a setting, to be used throughout the paper, 
\emph{time is counted backwards}, i.e., $d$ is the number of remaining
periods at which the project can be engaged. It turns out that
it is optimal to engage the latter project when it occupies state 
$i$ and $d$ periods remain  iff $\lambda^*(d, i) \geq \lambda$,
i.e., iff its current index 
is greater than or equal to the standard project's reward.
Such a result was first established for an undiscounted Bayesian Bernoulli
bandit  in \citet[Sect.\ 4]{bjk56}.
For an overview and extensions see \citet[Ch.\ 5]{berryFrist85}.

An economically insightul alternative representation of such an index is
\begin{equation}
\label{eq:orgirep}
\lambda^*(d, i) = \max_{1 \leq \tau \leq d} 
\frac{\displaystyle \Ex_{i}^{\tau}\left[\sum_{t=0}^{\tau-1}
   \beta^t  R\big(X(t)\big)\right]}{\displaystyle \Ex_{i}^{\tau}\left[\sum_{t=0}^{\tau-1} \beta^t\right]},
\end{equation}
where the right-hand-side is an
\emph{optimal-stopping problem}, with
$\tau$ denoting a 
\emph{stopping-time rule} for abandoning
the project provided it is engaged at least once starting at 
$i$ with $d$ remaining periods.
Thus,  $\lambda^*(d, i)$ is an \emph{average productivity} (AP) index, 
measuring the
maximum rate of expected  discounted reward that can be earned per unit of
 expected  
discounted time expended 
 by successively engaging the project no
more than $d$ times
starting at $i$.

Besides characterizing optimal policies for the finite-horizon 
two-armed bandit problem with one arm known, $\lambda^*(d, i)$ also serves
as a \emph{dynamic priority index} for engaging a project,  furnishing
a heuristic \emph{index policy} for the general FHMABP (\ref{eq:dsp}) that
engages at each time a project of currently highest index
value.
In the problem variant where \emph{at most} $K < M$ projects are to be
engaged at each time, the resulting index policy 
engages the project(s) with larger positive index values, if any, 
up to a maximum of
$K$.
Such a variant is particularly relevant when \emph{observation
  costs} or \emph{activity charges} are incorporated into the model. 
Thus, note that it follows immediately 
from (\ref{eq:orgirep}) that, if the aforementioned project model is modified to include 
a charge $\lambda$ to be incurred each time the project
is engaged, so the active reward is $\widetilde{R}(i) \triangleq R(i) -
\lambda$, the corresponding index becomes $\widetilde{\lambda}^*(d, i) = 
\lambda^*(d, i) - \lambda$.

The empirical performance of the index rule based
on $\lambda^*(d, i)$
for the case of two Bernoulli
projects with Beta priors is investigated in \cite{ginClayt99}, 
where it is shown to be very close to optimal for small time
horizons. 
In  \cite{carogall07}, such an index
rule is also considered, although the authors introduce and
use instead
 an approximation $\widehat{\lambda}^*(d, i)$ for $\lambda^*(d, i)$
 based on \emph{approximate DP}, which
is less costly to evaluate.

The index $\lambda^*(d, i)$ is monotone nondecreasing in the remaining time
$d$. Hence, for a project with bounded rewards it  has a
finite limit $\lambda^{*}(i)$ as $d \to \infty$. 
\cite{bellman56} showed that $\lambda^{*}(i)$
characterizes optimal
policies for the infinite-horizon two-armed Bernoulli bandit problem with one
arm known and $\beta < 1$. 
The resulting index rule was shown in \cite{gijo74} to be
optimal for the infinite-horizon multiarmed bandit
problem with one project engaged at each time, which has led
to $\lambda^*(i)$ being known as the \emph{Gittins index}. 

Although efficient algorithms to compute the Gittins index 
  of a  finite-state project are available, the currently lowest  \emph{time complexity} --- counting the number of \emph{arithmetic operations} (AOs) ---
for a general $n$-state project being $(2/3) n^3 + O(n^2)$ (as Gaussian elimination) for the algorithm given in
 \cite{nmijoc07}, 
the Gittins index $\lambda^{*}(i)$ of a
countably-infinite state project can only be approximated, using
  the
AP index $\lambda^*(d, i)$ for a large horizon
$d$.
\cite{wang97} shows the rate of convergence of $\lambda^*(d, i)$ to $\lambda^{*}(i)$ to be  linear for $\beta < 1$.

In contrast, 
scant research attention has been given to the efficient
 computation of the AP index  $\lambda^*(d, i)$ for a general project. 
The limited previous work, which we review in 
Section \ref{s:cfhi}, 
 typically focuses on specific models, and either uses 
DP to obtain approximate index values, or draws on the optimal-stopping
representation 
(\ref{eq:orgirep}) to obtain exact index values.
The latter approach, however,  has not yielded an 
exact index algorithm of general applicability with polynomial time complexity in both the horizon
$T$ and in the number of states $n$. 

A one-pass exact index
 algorithm with  $O(T^3 n^3)$ time complexity was presented in
\cite{nmcdc05}, using the \emph{adaptive-greedy} restless-bandit index algorithm 
introduced in \cite{nmaap01,nmmp02}.
 The term ``adaptive-greedy'' refers to the fact that such an
algorithm finds
a local maximizer for a certain vector at each step 
in an \emph{adaptive} fashion, meaning that
 such a vector is updated after the step.
A recursive version with an improved $O(T^2 n^3)$ time complexity, using $O(T^2
n^2)$ memory locations, was 
presented in 
\cite{nmsmct08}, by exploiting the structure of
the restless bandit formulation of a finite-horizon bandit.

This paper develops, simplifies, and extends such work, as it does not rely 
on restless bandit indexation and extends the algorithm's
scope to countably-infinite state projects. 
For a project with limited
transitions per state, the
 time complexity is reduced to $O(T^2 n^2)$, and in such a case the
 algorithm is further extended to countably-infinite state projects with an $O(T^6)$ time complexity and an $O(T^5)$ memory complexity.
A computational study demonstrates the index algorithm's practical
tractability for moderate-size instances.

For comparison, 
the paper includes in Section \ref{s:caaic} an assessment of the
complexity of approximate index computation via the conventional 
\emph{calibration method}, which solves by DP a collection
of finite-horizon optimal stopping problems
 at a finite grid of $\lambda$-values (from which the approximate
 index values are taken). 
The complexity of such a method grows linearly in the grid size
 $L$, being $O(L T n^2)$ time and $O(L T n)$
space.
Hence, if $L$ is taken equal to the number $T n$ of index values to be evaluated, 
one obtains time and space complexities of $O(T^2 n^3)$ and 
$O(T^2 n^2)$, as in the exact algorithm proposed herein.

The remainder of the paper is organized as follows.
Section \ref{s:cfhi} reviews previous approaches to the finite-horizon
AP index' computation.
Section \ref{s:fic} develops the recursive 
index algorithm for finite-state projects. 
Section \ref{s:icism} extends the algorithm's scope to
countably-infinite 
state projects.
Section \ref{s:riraga} presents  an efficient \emph{block implementation} (see \citet{dongeijk00}) of the algorithm, which 
is necessary for making it useful in practice.
Section \ref{s:ce} reports the results of a computational study.
Section \ref{s:concl} concludes.

Ancillary material is available in an online supplement.

\section{Previous approaches to the AP index computation}
\label{s:cfhi}
Two approaches have been proposed, as discussed below.
\subsection{The calibration method for approximate index computation}
\label{s:caaic} 
The first, known as the \emph{calibration method}, uses DP to obtain approximate index values,
adapting to the finite-horizon setting the method in 
\citet[Sec.\ 8]{gi79} for approximate Gittins index calculation, 
as outlined in \citet[Ch.\
5]{berryFrist85}.
 Denoting by 
$v_d^*(i; \lambda)$ the optimal value of the two-armed bandit problem
with one arm known, where the unknown project starts at $i$ with
$d$ remaining periods, $\lambda^*(d, i)$ is the 
smallest root in $\lambda$ of 
\begin{equation}
\label{eq:lambdadefeq}
v_d^*(i; \lambda) = \lambda h_d,
\end{equation}
where $h_d \triangleq (1-\beta^d)/(1-\beta)$ if $\beta < 1$ and 
$h_d  \triangleq d$ if $\beta = 1$.
Note that, for a fixed $\lambda$, 
 $v_d^*(i; \lambda)$ is recursively characterized by the DP equations
\begin{equation}
\label{eq:vddpe}
v_d^*(i; \lambda) = 
\begin{cases}
\max \big\{\lambda h_d, R(i)  + 
  \beta \sum_{j \in \mathbb{X}} p(i, j) v_{d-1}^*(j; \lambda)\big\},
& \quad d \geq 2 \\
\max \big\{\lambda, R(i)\big\}, & \quad d = 1,
\end{cases}
\end{equation}
which use the result, first proven in \citet[Lemma 4.1]{bjk56},
 that, if the standard project is optimal at any stage, then it is also
 optimal thereafter.

The calibration method solves such DP equations for a
grid  of increasing $\lambda$-values $\{\lambda_l\colon 1 \leq l \leq L\}$,  with $\lambda_1 = \min_i R(i)$ and $\lambda_L = \max_i R(i)$, which gives index approximations  $\widehat{\lambda}^*(d, i)$ 
with the desired degree of accuracy. 
Table \ref{tab:cm} shows an efficient \emph{block implementation} (see \citet{dongeijk00}) of the calibration method, 
where $\mathbf{V}_d^*$ is the $n \times L$ matrix $[\mathbf{v}_d^*(\lambda_l)]_{1 \leq l \leq L}$, with 
$\mathbf{v}_d^*(\lambda_l)  
 = \big(v_d^*(i; \lambda_l)\big)_{i \in \mathbb{X}}$, and
$\mathbf{1}$ is an $n$-vector of ones. 
Note that the ``max" shown in Table \ref{tab:cm} are to be read componentwise.
As discussed in Section \ref{s:riraga}, block implementations achieve economies of scale  in computation by rearranging bottleneck calculations as 
operations on large data blocks. In Table \ref{tab:cm}, this is achieved with the matrix-update $\mathbf{V}_d^* :=  \beta \mathbf{P} \mathbf{V}_{d-1}^*$ . 
The required minimization can be carried out using bisection search.

The following result assesses both the time (AOs) and
the memory complexities (the latter measuring intermediate floating-point
storage locations, i.e., excluding input and
output)
of the calibration method. 
Since it appears reasonable to deploy such an approach
using a grid contaning a number $L$ of 
$\lambda$-values that is at least as large as the number of
 index values, i.e., $L \geq T n$, Proposition 
\ref{pro:caaoc}(c) estimates the complexity in the case $L = T n$.

\begin{table}[htb]
\caption{Block implementation of the calibration method.}
\begin{center}
\fbox{%
\begin{minipage}{4in}
\textbf{ALGORITHM} $\mathrm{BlockCAL}$: \\
\textbf{Input:} $\{\lambda_l\colon 1 \leq l \leq L\}$ \\
\textbf{Output:} $\{\widehat{\lambda}^*(d, i)\colon 1 \leq d \leq T, i \in \mathbb{X}\}$ 
\begin{tabbing}
$\widehat{\lambda}^*(1, i) := R(i), \, i \in \mathbb{X}$; \, $\mathbf{v}_1^*(\lambda_l) :=  \max \big\{\lambda_l \mathbf{1}, \mathbf{R}\big\}, \, l = 1, \ldots, L$ \\
\textbf{for} \= $d := 2$ \textbf{to}  $T$ \textbf{do} \\
\> $\mathbf{V}_d^* := \beta \mathbf{P} \mathbf{V}_{d-1}^*$ \{ note:
$\mathbf{V}_d^* = [\mathbf{v}_d^*(\lambda_l)]_{1 \leq l \leq L}$ \} \\
\> $\mathbf{v}_d^*(\lambda_l) :=  \max \big\{h_d \lambda_l \mathbf{1}, \mathbf{R} + \mathbf{v}_{d}^*(\lambda_l)\big\}, \, l = 1, \ldots, L$ \\
\> $\widehat{\lambda}^*(d, i) := \min \big\{\lambda_l\colon v_d^*(i; \lambda_l) = \lambda_l h_d, 1 \leq l \leq L\big\}, \, i \in \mathbb{X}$ \\
\textbf{end} \{ for \}
\end{tabbing}
\end{minipage}}
\end{center}
\label{tab:cm}
\end{table}

\begin{proposition}
\label{pro:caaoc} For an $n$-state project with horizon $T:$
\begin{itemize}
\item[\textup{(a)}]  For a fixed $\lambda$,
the $\mathbf{v}_d^*(\lambda)$ for 
$1 \leq d \leq T$  can be computed in
$2 (T-1) \big[n (n+1) + 1\big] + n = O(T n^2)$ time and $O(T n)$
space.
\item[\textup{(b)}] Using a size-$L$ grid, the calibration method
 uses $2 (T-1) L \big[n (n+1) + 1\big] + L n =
 O(L T n^2)$ time and
 $O(L T n)$ space.
\item[\textup{(c)}] If 
$L = T n$, then the  calibration method uses $O(T^2 n^3)$ time and
 $O(T^2 n^2)$ space.
\end{itemize}
\end{proposition}
\begin{proof}
(a)
Computing $v_1^*(\cdot; \lambda)$ involves the $n$ subtractions
required to obtain $\max \big\{\lambda, R(i)\big\}$ for every $i$.
For $d \geq 2$, in order to compute $v_d^*(i; \lambda)$ given $v_{d-1}^*(\cdot; \lambda)$ and 
$\tilde{h}_{d-1} = \lambda h_{d-1}$, one first computes $\sum_{j} p(i, j) v_{d-1}^*(j;
\lambda)$, which takes $2 n - 1$ AOs. Multiplying the result  by 
$\beta$, adding it to $R(i)$, and then subtracting $\tilde{h}_d$ to
determine the ``max'' takes 3 AOs. Repeating for every $i$ takes
$2 n (n+1)$ AOs. The $\tilde{h}_d$ is computed from
$\tilde{h}_d = \lambda + \beta \tilde{h}_{d-1}$, which takes 2
additional AOs.
Repeating for $d = 2, \ldots, T$ gives the $2 (T-1) \big[n (n+1) +
1\big]$ term. Further, storage of the $v_d^*(\cdot; \lambda)$ uses 
$T n$ floating-point memory locations.

Parts (b) and (c) follow immediately from parts (a) and (b), respectively.
\end{proof}

Note further that the calibration method is immediately parallelizable, as the computations for nonoverlapping ranges of $\lambda$-values can be split
among different processors. Hence, its time complexity scales linearly with the number of processors.

\subsection{The direct method for exact  index computation}
\label{s:eica}
In contrast to the calibration method, which  is presently the preferred approach, 
the  direct method computes exact index values and is relatively 
unexplored.  It was  introduced in \citet[Sec.\
4]{bjk56}, and has been extended in \citet[Ch.\ 5]{berryFrist85} to more
 general discount sequences. 
The direct method draws on
the representation in (\ref{eq:orgirep}), calling for the solution of the corresponding
optimal stopping problems.

In \citet[Sec.\ 7]{gi79}, such a method is deployed to  compute the
index $\lambda^*(T, i)$ for a Bernoulli bandit with Beta priors,
with the goal  of approximating its Gittins index $\lambda^*(i)$.
Two key results, which  reduce
 the complexity of the 
optimal stopping problems of concern, and give a recursive index
computation,
are Corollaries 1 and 2 in that paper.

\begin{proposition}[Corollary 1 in \cite{gi79}]
\label{pro:cor1gi79}
An optimal stopping time for 
 \textup{(\ref{eq:orgirep})} is  
\begin{equation}
\label{eq:gi79cor1}
\tau_d^* \triangleq \min \Big\{d, \min \big\{t\colon 1 \leq t \leq d-1, 
\lambda^*\big(d-t, X(t)\big) < \lambda^*(d, i)\big\}\Big\}.
\end{equation}
\end{proposition} 

\begin{proposition}[Corollary 2 in \cite{gi79}]
\label{pro:cor2gi79}
In order to solve optimal-stopping
problem \textup{(\ref{eq:orgirep})}, it suffices to
 consider stopping times of the form
\begin{equation}
\label{eq:taulambda}
\tau_d(\lambda) \triangleq \min \Big\{d, \min \big\{t\colon 1 \leq t \leq d-1, 
\lambda^*\big(d-t, X(t)\big) < \lambda\big\}\Big\}, \quad \lambda \in \mathbb{R},
\end{equation}
where the ``\textup{min}'' of an empty set is taken to be $\infty$.
\end{proposition} 

The direct  method outlined in \citet[Sec.\ 7]{gi79} 
draws on Proposition \ref{pro:cor1gi79}, 
 reducing optimal-stopping problem (\ref{eq:orgirep}) to the 
one-dimensional continuous optimization problem
\begin{equation}
\label{eq:lambdadtaulamb}
\lambda^*(d, i) = \max_{\lambda \in \mathbb{R}} 
\frac{\displaystyle \Ex_{i}^{\tau_d(\lambda)}\left[\sum_{t=0}^{\tau_d(\lambda)-1}
    \beta^t R\big(X(t)\big)\right]}{\displaystyle
  \Ex_{i}^{\tau_d(\lambda)}\left[\sum_{t=0}^{\tau_d(\lambda)-1} \beta^t\right]},
\end{equation}
which corresponds to formula (11) in that paper.

However, in \cite{gi79} it is not discussed how to solve exactly the 
right-hand-side optimization problem over $\lambda \in \mathbb{R}$ 
in (\ref{eq:lambdadtaulamb}), although in the review of such an
approach in \citet[p.\ 140]{gi89} it appears to be suggested to use
a grid of $\lambda$-values and interpolation for such a purpose, 
which would render the method approximate rather than exact. 

\citet[Sec.\ IV]{vawabu85} gives a 
 Gittins-index algorithm that exploits the corresponding
result in Proposition \ref{pro:cor1gi79} for the Gittins index, as
stated in the lemma on p.\ 154 of \cite{gi79}. The Varaiya et al.\  (1985) algorithm
avoids solving the corresponding continuous optimization problem
(\ref{eq:lambdadtaulamb}), as it only involves discrete maximizations.

In fact, Proposition \ref{pro:cor1gi79} ensures that the continuous
optimization problem in the right-hand-side of 
(\ref{eq:lambdadtaulamb}) can be reduced to the  
discrete optimization problem
\begin{equation}
\label{eq:doplambdadtaulamb}
\lambda^*(d, i) = \max_{\lambda \in \{\lambda^*(s, j)\colon 1 \leq s
  \leq d-1, j \in \mathbb{X}\}} 
\frac{\displaystyle \Ex_{i}^{\tau_d(\lambda)}\left[\sum_{t=0}^{\tau_d(\lambda)-1}
     \beta^t R\big(X(t)\big)\right]}{\displaystyle
  \Ex_{i}^{\tau_d(\lambda)}\left[\sum_{t=0}^{\tau_d(\lambda)-1} \beta^t\right]}.
\end{equation}
Yet, this observation does not directly yield  an adaptive-greedy
algorithm for the AP index $\lambda^*(d, i)$ 
analogous to that of Varaiya et al.\ for the Gittins index
$\lambda^*(i)$.

\section{Recursive index computation}
\label{s:fic}
This section develops the recursive adaptive-greedy index algorithm
that is the main contribution of this paper, 
for a project with a finite number $n$ of states.
\subsection{Reduction  to a modified Gittins index}
\label{s:icpclagia}
Let us first  prepare the ground for computing
a project's finite-horizon AP index $\lambda^*(d, i)$,
by showing that such an index can be reduced to a modified Gittins index, which allows 
use of the adaptive-greedy algorithm available for the latter index to
compute the former.

Consider an auxiliary \emph{infinite-horizon} project, whose state 
$Y(t)$ evolves over time periods $t \geq 0$ through the state space 
$\mathbb{Y}_T \triangleq \mathbb{Y}_T^{\{0, 1\}} \cup
\mathbb{Y}^{\{0\}}$, where $\mathbb{Y}_T^{\{0, 1\}}
\triangleq \{1, \ldots, T\} \times \mathbb{X}$ is the set of 
\emph{controllable states} where both actions (active and passive) are
allowed, and 
$\mathbb{Y}^{\{0\}} \triangleq \{(0, \Omega)\}$ is the (singleton) set
of \emph{uncontrollable states} where the passive action \emph{must}
be taken, with $(0, \Omega)$ denoting a terminal absorbing state.
Under the active action $a(t) = 1$, the project's
transition
probabilities are 
$p^1\big((d, i), (d-1, j)\big) \triangleq p(i, j)$ for $2 \leq d \leq T$ and 
$p^1\big((1, i), (0, \Omega)\big) \triangleq 1$, while its
rewards
are 
$R^1(d, i) \triangleq R(i)$.
Under the passive action $a(t) = 0$, 
the project remains frozen,  its
transition
probabilities being 
$p^0\big((d, i), (d, i)\big) \equiv 1$, while its immediate rewards
are $R^0(d, i) \equiv 0$. Further,  all other transition probabilities are zero.

The idea of forcing a project to be passive in certain states,
termed uncontrollabe,  
 was introduced in 
\cite{nmmp02} in the setting of restless bandits,
where  a project's index is only defined for its
controllable states.
This
is relevant in the present setting, as shown next.
Suppose we allow the active action to be taken at the absorbing state $(0, \Omega)$  in the auxiliary infinite-horizon project described above, with 
 the same dynamics and rewards as the passive action.
Let $G(d, i)$ be the corresponding \emph{conventional Gittins index}, 
which is defined   by 
\begin{equation}
\label{eq:gsic}
G(d, i) \triangleq \max_{\tau \geq 1} \frac{\displaystyle \Ex_{(d,
    i)}^{\tau}\left[\sum_{t=0}^{\tau-1} \beta^t R^1\big(Y(t)\big)\right]}{\displaystyle \Ex_{(d,
    i)}^{\tau}\left[\sum_{t=0}^{\tau-1}  \beta^t\right]}, \quad (d, i)
\in \mathbb{Y}_T,
\end{equation}
where $\tau \geq 1$ is a stopping time under which the project is
engaged at least once starting at $(d, i)$.
Now, let $G'{(d, i)}$ be the \emph{modified Gittins index}
 of the project with state $Y(t)$
and uncontrollable state $(0, \Omega)$, which is defined only for its
controllable states $(d, i) \in \mathbb{Y}_T^{\{0, 1\}}$ by 
\begin{equation}
\label{eq:gsi}
G'{(d, i)} \triangleq \max_{\tau \geq 1\colon a(t) = 0 \text{ if } Y(t) = 
 (0, \Omega)} \frac{\displaystyle \Ex_{(d,
   i)}^{\tau}\left[\sum_{t=0}^{\tau-1} \beta^t R^1\big(Y(t)\big)
   \right]}{\displaystyle \Ex_{(d,
   i)}^{\tau}\left[\sum_{t=0}^{\tau-1}  \beta^t\right]} = 
\max_{1 \leq \tau \leq d} \frac{\displaystyle \Ex_{(d,
   i)}^{\tau}\left[\sum_{t=0}^{\tau-1} \beta^t R^1\big(Y(t)\big)
   \right]}{\displaystyle \Ex_{(d,
   i)}^{\tau}\left[\sum_{t=0}^{\tau-1}  \beta^t\right]}.
\end{equation}

Note that, unlike (\ref{eq:gsic}), optimal-stopping problem  (\ref{eq:gsi}) only considers
 stopping times $\tau \geq 1$  that  idle the project  at $(0, \Omega)$, i.e.,  $\tau \leq
d$.
Such a distinction between the conventional and the modified Gittins index
is significant, since 
in some cases the two may differ.
Thus, e.g., for a state $(1, i)$ with $R(i) < 0$, 
$G'{(1, i)} = R(i) < G{(1, i)} = (1-\beta) R(i)$.

The interest of introducing such an auxiliary project and its modified
Gittins index $G'{(d, i)}$ is that the latter is precisely the
finite-horizon AP index $\lambda^*{(d, i)}$ of concern here.

\begin{proposition}
\label{pro:rgi}
$\lambda^*{(d, i)} = G'{(d, i)}$ for $(d, i) \in  \mathbb{Y}_T^{\{0, 1\}}$.
\end{proposition}
\begin{proof}
The result follows by noting that, under a stopping time 
$1 \leq \tau \leq d$
for the original project's state 
process $X(t)$ starting at $i$ with $d$ remaining periods, 
the process defined by
 $Y(t) \triangleq \big(d-t, X(t)\big)$ for $t = 0, \ldots, \tau-1$ has
the same (active) dynamics and rewards as that used to define $G'{(d, i)}$
above, and therefore
\[
\lambda^*{(d, i)} \triangleq 
\max_{1 \leq \tau \leq d} \frac{\displaystyle
  \Ex_{i}^{\tau}\left[\sum_{t=0}^{\tau-1} \beta^t R\big({X(t)}\big)
    \right]}{\displaystyle
  \Ex_{i}^{\tau}\left[\sum_{t=0}^{\tau-1}  \beta^t\right]} = 
\max_{1 \leq \tau \leq d} \frac{\displaystyle \Ex_{(d,
   i)}^{\tau}\left[\sum_{t=0}^{\tau-1} \beta^t R^1\big({Y(t)}\big)
   \right]}{\displaystyle \Ex_{(d,
   i)}^{\tau}\left[\sum_{t=0}^{\tau-1}  \beta^t\right]} = 
G'{(d, i)}.
\]
\end{proof}

\subsection{Adaptive-greedy index algorithm}
\label{s:agia}
The representation in Proposition \ref{pro:rgi} of the finite-horizon
index $\lambda^*{(d, i)}$
as the modified Gittins index of
 an infinite-horizon project allows us to use
available algorithms for the latter type of index to compute 
the former.
Note, however, that the classic adaptive-greedy Gittins-index algorithm of 
\cite{vawabu85} should not be directly used, since it computes the 
conventional Gittins index $G{(d, i)}$ which, as argued above, can
differ from the modified Gittins index $G'{(d, i)} =\lambda^*{(d, i)}$.
Also, such an algorithm does not exploit special structure. 

We will use instead an extension of such an algorithm
introduced in \cite{nmaap01} for restless
bandits, and further extended in \cite{nmmp02} to a wider setting. A key feature of
 such 
an algorithm  is that it exploits special structure to reduce the computational
burden.
Rather than describing the algorithm in full generality,
for which the reader is referred to the aforementioned papers, 
we present it next as it applies to the model of concern. 

To prepare the ground, we start by 
introducing two measures to evaluate a stopping-time rule $0 \leq \tau
\leq d$ for the project (refer to the discussion in
Section \ref{s:icpclagia}): a
\emph{reward measure}
\[
f_d^{\tau}{(i)} \triangleq \Ex_{(d, i)}^{\tau}\left[\sum_{t=0}^{\tau-1} 
  \beta^t  R^{a(t)}\big({Y(t)}\big)\right] = 
\Ex_{i}^{\tau}\left[\sum_{t=0}^{\tau-1} 
   \beta^t R\big({X(t)}\big)\right],
\]
giving the expected total discounted reward earned
starting at $i$ with $d$ remaining periods; and the
\emph{work measure}
\[
g_d^{\tau}{(i)} \triangleq \Ex_{(d, i)}^{\tau}\left[\sum_{t=0}^{\tau-1} \beta^t\right],
\]
giving the corresponding
 expected total discounted time that the project is active. 

Due to the optimality of deterministic  Markov policies 
for finite-state and -action MDPs, it suffices to consider stopping
times $\tau$ given by a \emph{continuation} (or \emph{active}) \emph{set} 
$A \subseteq \mathbb{Y}_T^{\{0, 1\}}$, consisting of those
controllable states at which the project is 
active under $\tau$.
We will find it convenient to represent each such continuation set
in a more explicit fashion, writing
\begin{equation}
\label{eq:adec}
A = (A_1, \ldots, A_T) \triangleq \{1\} \times A_1 \cup \cdots \cup \{T\}
  \times A_T,
\end{equation}
where $A_d \in \mathbb{X}$ is the continuation set when $d$ periods
remain. 
Thus, the stopping rule having continuation set $A$ 
engages the original project in state $i$ when $d$ periods remain  (or the modified project in state $(d, i)$) iff 
$i \in A_d$.
We will  write $f^A_d{(i)}$ and $g^A_d{(i)}$  to denote the
   reward measures and work measures, respectively, under such a stopping rule.

Further, we will use the \emph{modified reward} and \emph{work measures}
defined by 
\begin{equation}
\label{eq:rawasi}
r^A_d{(i)} \triangleq f^{A \cup \{(d, i)\}}_d{(i)}, \quad \text{and} \quad
w^A_d{(i)} \triangleq g^{A \cup \{(d, i)\}}_d{(i)},
\end{equation}
respectively, along with the \emph{productivity rate measure} defined by 
\begin{equation}
\label{eq:lambdaasi}
\lambda^A_d{(i)} \triangleq \frac{r^A_d{(i)}}{w^A_d{(i)}}.
\end{equation}

Now, the classic result referred to in Section \ref{s:intro}, whereby if it is
optimal  to stop the project when $d$ periods
remain, then it is also optimal to do so when less
periods remain, allows us to restrict the
continuation sets that need be considered to those consistent with
such a property, which constitute the \emph{continuation-set
  family} 
\begin{equation}
\label{eq:fhat}
\mathcal{F}_T \triangleq \big\{A = (A_1, \ldots, A_T)\colon A_1  \subseteq \cdots \subseteq
A_T \subseteq \mathbb{X}\big\}.
\end{equation}

We are now ready to present the \emph{adaptive-greedy index
  algorithm} $\mathrm{AG}({\mathcal{F}_T})$, which is shown in Table
\ref{tab:agaf}.
Such an algorithm  builds up
in $T n$ steps (note that $T n = T |\mathbb{X}| = |\mathbb{Y}_T^{\{0, 1\}}|$ is
the number of  controllable states)
an increasing nested chain of \emph{adjacent} continuation sets (i.e., differing by one
state) $A^0 = \emptyset \subset
A^1 \subset \ldots \subset A^{T n} = \mathbb{Y}_T^{\{0, 1\}}$ in
$\mathcal{F}_T$ connecting the empty set to 
the full controllable state space, proceeding at each step in a greedy fashion.
Thus, once the  continuation set $A^{k-1} \in \mathcal{F}_T$ has been constructed, the next continuation
set $A^k$ is obtained by adding to $A^{k-1}$ a controllable
state $(s^k, i^k) \in \mathbb{Y}_T^{\{0, 1\}} \setminus A^{k-1}$ that
maximizes 
the productivity rate 
$\lambda_s^{A^{k-1}}{(i)}$ over augmented states $(s, i) \in \mathbb{Y}_T^{\{0, 1\}} \setminus A^{k-1}$ for which the next active
set remains in $\mathcal{F}_T$, i.e., with $A^k = A^{k-1}
\cup \{(s, i)\} \in \mathcal{F}_T$. Ties are broken arbitrarily.

Note that in Table \ref{tab:agaf} we write, for notational convenience,
$\lambda_s^{A^{k-1}}{(i)}$ as
$\lambda_s^{k-1}{(i)}$.
The algorithm's output consists of an augmented-state sequence $(s^k, i^k)$
spanning $\mathbb{Y}_T^{\{0, 1\}}$, along with a corresponding  nonincreasing sequence of index values
  $\lambda^*(s^k, i^k)$.

\begin{table}[htb]
\caption{Adaptive-greedy index algorithm 
$\mathrm{AG}({\mathcal{F}_T})$.}
\begin{center}
\fbox{%
\begin{minipage}{4in}
\textbf{ALGORITHM} $\mathrm{AG}({\mathcal{F}_T})$: \\
\textbf{Output:} $\{(s^k, i^k), \lambda^*(s^k, i^k)\colon 1 \leq k \leq T n\}$ 
\begin{tabbing}
$A^0 := \emptyset$ \\
\textbf{for} \= $k := 1$ \textbf{to}  $T n$ \textbf{do} \\
\> \textbf{pick}  
 $(s^*, i^*) \in \argmax
      \big\{\lambda^{{k-1}}_{s}(i)\colon
                 (s, i) \in \mathbb{Y}_T^{\{0, 1\}} \setminus A^{k-1}, 
     A^{k-1} \cup \{(s, i)\} \in \mathcal{F}_T\big\}$ \\
 \>  $\lambda^*(s^*, i^*) := 
 \lambda^{{k-1}}_{s^*}(i^*)$;  \, $A^{k} := A^{k-1} \cup \{(s^*, i^*)\}$; \, $(s^k, i^k) := (s^*, i^*)$ \\
\textbf{end} \{ for \}
\end{tabbing}
\end{minipage}}
\end{center}
\label{tab:agaf}
\end{table}

We next use the definition of 
$\mathcal{F}_T$ in (\ref{eq:fhat}) to obtain the more explicit
 reformulation of  algorithm $\mathrm{AG}({\mathcal{F}_T})$ shown in Table \ref{tab:expalg},  
 breaking down the choice
 at step $k$ of the augmented state 
to be added to the current continuation set $A^{k-1} = (A_1^{k-1},  \ldots, A_T^{k-1})$.

\begin{table}[htb]
\caption{More explicit reformulation of index algorithm
  $\mathrm{AG}({\mathcal{F}_T})$.}
\begin{center}
\fbox{%
\begin{minipage}{4in}
\textbf{ALGORITHM} $\mathrm{AG}({\mathcal{F}_T})$ \\
\textbf{Output:}
$\{(s^k, i^k), \lambda^*(s^k, i^k)\colon 1 \leq k \leq T n\}$
\begin{tabbing}
$A_1^0 := \cdots := A_T^0 := \emptyset$; $A_{T+1}^0 :=
\mathbb{X}$ \\
\textbf{for} \= $k := 1$ to $T n$
 \textbf{do} \{ note: $A^{k-1} = (A_1^{k-1}, \ldots,
 A_T^{k-1})$ \} \\
\> \textbf{pick} $(s^*, i^*)  \in \argmax \big\{\lambda^{k-1}_s{(i)}\colon 1 \leq s \leq T, A_s^{k-1} \subset A_{s+1}^{k-1}, i \in
A_{s+1}^{k-1} \setminus A_s^{k-1}\big\}$ \\
\>  $\lambda^*(s^*, i^*) := 
 \lambda^{{k-1}}_{s^*}(i^*)$;  $A_{s^*}^{k} := A_{s^*}^{k-1} \cup \{i^*\}$; 
  $A_{s}^{k} := A_{s}^{k-1}, s \neq s^*$; $(s^k, i^k) := (s^*, i^*)$ \\
\textbf{end }  \{ for \} 
\end{tabbing}
\end{minipage}}
\end{center}
\label{tab:expalg}
\end{table}

\subsection{Reward and work measure recursions}
\label{s:rmrm}
The above algorithms do not 
specify how to compute the required modified reward and work
measures  (see (\ref{eq:lambdaasi})) to calculate
 productivity rates $\lambda^{k-1}_d{(i)}$.
This section presents recursions that will be used for
such a purpose in the next section.

Let $A = (A_1, \ldots, A_T)$ be a 
continuation set in  $\mathcal{F}_T$.
Note first that, from the stopping-rule interpretation of  $A$,  it is clear
that the reward measure $f^A_d{(i)}$ does not depend on $A_{s}$ for
$s > d$, 
which allows us to write 
$f^A_d{(i)}$ as $f_d^{(A_1, \ldots, A_d)}{(i)}$ for $d
\geq 1$, while $f_0^A{(\Omega)} \equiv 0$.
Further, the definition of modified reward measure $r_d^A{(i)}$
in (\ref{eq:rawasi}) ensures that it does not depend on $A_{s}$ for
$s \geq d$, 
which allows us to write $r_d^A{(i)}$ as $r_d^{(A_1, \ldots,
  A_{d-1})}{(i)}$ for $d \geq 2$, while $r_1^A{(i)} = R(i)$.

In the following result,
 part (a) shows how to recursively evaluate modified reward
measures $r_d^A{(i)}$ for a fixed $A$. The remaining parts show how to 
evaluate the modified reward measure for an augmented continuation
set, $r_d^{A \cup \{(d, i^*)\}}{(i)}$, based on knowledge of the
$r_d^A{(i)}$.
While the base continuation set is written as $A$, the
reduction discussed in the previous paragraph should be taken into 
account, as it plays a key role in the proof. 
Also, note that $A \cup \{(d, i^*)\} = (A_1, \ldots, A_d \cup \{i^*\}, \ldots, A_T)$.
\begin{lemma}
\label{lma:mrmr}
For $1 \leq d \leq T$, $i \in \mathbb{X}$\textup{:} 
\begin{itemize}
\item[\textup{(a)}] 
$r_d^{A}{(i)} = 
\begin{cases}
\displaystyle
R(i)
 + \beta \sum_{j \in A_{d-1}} p({i, j}) r_{d-1}^{A}{(j)} & \text{ if
 }  2 \leq d \leq T \\
R(i) & \text{ if } d = 1.
\end{cases}
$
\item[\textup{(b)}] 
$r_d^{A \cup \{(d, i^*)\}}{(i)} = 
r_d^{A}{(i)}$, for $i^* \in \mathbb{X} \setminus A_d$.
\item[\textup{(c)}] 
$r_d^{A \cup \{(d-1, i^*)\}}{(i)}  = 
r_d^{A}{(i)} + 
\beta p({i, i^*}) r_{d-1}^{A}(i^*)$, 
for $i^* \in A_d \setminus A_{d-1}$.
\item[\textup{(d)}] 
$\displaystyle r_d^{A \cup \{(s, i^*)\}}{(i)} = 
R(i) + 
\beta \sum_{j \in A_{d-1}} p({i, j}) r_{d-1}^{A
  \cup \{(s, i^*)\}}(j)$, 
for $1 \leq s \leq d-2$, $i^* \in A_{s+1} \setminus A_{s}$.
\end{itemize}
\end{lemma}
\begin{proof}
(a) This part follows immediately from the definition of 
$r_d^A{(i)}$ in (\ref{eq:rawasi}) and the project's dynamics and
rewards under the stopping rule induced by continuation set $A$.

(b) The result follows from 
\[
r_d^{A \cup \{(d, i^*)\}}{(i)}  =
r_{d}^{(A_1, \ldots, A_{d-1})}(i) = r_{d}^A(i).
\]

(c) The result, which draws on part (a),  follows from
\begin{align*}
r_d^{A \cup \{(d-1, i^*)\}}{(i)}  & = 
r_{d}^{(A_1, \ldots, A_{d-1} \cup \{i^*\})}(i)  = 
R(i) + \beta \sum_{j \in A_{d-1}} p({i, j}) r_{d-1}^{(A_1, \ldots, A_{d-1}
  \cup \{i^*\})}(j) \\
& \qquad + \beta p({i, i^*}) r_{d-1}^{(A_1, \ldots, A_{d-1}
  \cup \{i^*\})}(i^*) \\
& = R(i) + \beta \sum_{j \in A_{d-1}} p({i, j}) r_{d-1}^{(A_1,
  \ldots, A_{d-2})}(j) + \beta p({i, i^*}) r_{d-1}^{(A_1, \ldots,
  A_{d-2})}(i^*) \\
& = 
r_d^{A}{(i)} + 
\beta p({i, i^*}) r_{d-1}^A(i^*).
\end{align*}

(d) The result, which also draws on part (a),  follows from
\begin{align*}
r_{d}^{A \cup \{(d, i^*)\}}(i)  & = 
r_{d}^{(A_1, \ldots, A_{s} \cup \{i^*\}, \ldots, A_{d-1})}(i) = R(i)
+ \beta \sum_{j \in A_{d-1}} p({i, j}) r_{d-1}^{(A_1,
  \ldots, A_{s} \cup \{i^*\}, \ldots, A_{d-2})}(j) \\
& = R(i) + \beta \sum_{j \in A_{d-1}} p({i, j}) r_{d-1}^{A \cup \{(s,
  i^*)\}}(j).
\end{align*}
\end{proof}

The following result is the counterpart of Lemma \ref{lma:mrmr}
for modified work measures $w^{A}_d(i)$. It follows immediately from the above by 
taking $R(i) \equiv 1$, and hence we omit its proof.

\begin{lemma}
\label{lma:mwmr}
For $1 \leq d \leq T$, $i \in \mathbb{X}$\textup{:} 
\begin{itemize}
\item[\textup{(a)}] 
$w^{A}_{d}(i) = 
\begin{cases}
\displaystyle
1
 + \beta \sum_{j \in A_{d-1}} p({i, j}) w^{A}_{d-1}(j) & \text{ if }  d \geq 2 \\
1 & \text{ if } d = 1.
\end{cases}
$
\item[\textup{(b)}] 
$w_d^{A \cup \{(d, i^*)\}}{(i)}  = 
w_d^{A}{(i)}$, for $i^* \in \mathbb{X} \setminus A_d$,
\item[\textup{(c)}] 
$w_d^{A \cup \{(d-1, i^*)\}}{(i)}  = 
w_d^{A}{(i)} + 
\beta p({i, i^*}) w_{d-1}^{A}(i^*)$, 
for $i^* \in A_d \setminus A_{d-1}$;
\item[\textup{(d)}] 
$\displaystyle w_d^{A \cup \{(s, i^*)\}}{(i)}  = 
1 + 
\beta \sum_{j \in A_{d-1}} p({i, j}) w_{d-1}^{A
  \cup \{(s, i^*)\}}(j)$, 
for $1 \leq s \leq d-2$, $i^* \in A_{s+1} \setminus A_{s}.$
\end{itemize}
\end{lemma}

\subsection{A $T$-stage  $O(T^2 n^3)$ recursive  adaptive-greedy index algorithm}
\label{s:ria}
This section draws on the above results to reformulate the one-pass adaptive-greedy
 index algorithm
$\mathrm{AG}({\mathcal{F}_{T}})$ in Table \ref{tab:expalg} into a $T$-stage \emph{recursive adaptive-greedy} (RAG) 
algorithm.

\begin{table}[htb]
\caption{Stage $2 \leq d \leq T$ of the RAG index algorithm.}
\begin{center}
\fbox{%
\begin{minipage}{6in}
\textbf{ALGORITHM} $\mathrm{RAG}_{d}$ \\
\textbf{Input:}
$\big\{\lambda^*(s, i)\colon 1 \leq s \leq d-1, i \in \mathbb{X}\big\}, \big\{(s_{d-1}^l, i_{d-1}^l), \mathbf{w}^{l-1}_{d-1}, \mathbf{r}^{l-1}_{d-1}\colon 1 \leq l \leq L_{d-1}\}$
 \\
\textbf{Output:}
$\big\{\lambda^*(s, i)\colon 1 \leq s \leq d, i \in \mathbb{X}\big\}, \big\{(s_{d}^k, i_{d}^k), \mathbf{w}^{k-1}_d, \mathbf{r}^{k-1}_d\colon 1 \leq k \leq L_d\big\}$
\begin{tabbing}
$k := 1$; \, $l := 1$; \, $A_{d-1} := A_d := \emptyset$; \,   
$k_d := 0$; \, 
$\begin{bmatrix} \mathbf{w}^{0}_{d} & \mathbf{r}^{0}_d\end{bmatrix} :=  \begin{bmatrix} \mathbf{1} &
  \mathbf{R} \end{bmatrix}$ \\
\textbf{repeat} \=    \{ note: $A^{k-1} = (A_1^{k-1}, \ldots, A_d^{k-1}) = (A_1, \ldots, A_d)$ \} \\
\> $\lambda^{k-1}_{d}(i) := r^{k-1}_{d}(i)/w^{k-1}_{d}(i),  i \in \mathbb{X} \setminus A_d$;  \textbf{pick} 
 $\displaystyle i^{*} \in \argmax
\big\{\lambda^{k-1}_{d}(i)\colon i \in \mathbb{X} \setminus A_d\big\}$ \\
\> \textbf{if} \= $\lambda^{k-1}_{d}(i^*) \geq \lambda^*(s_{d-1}^l, i_{d-1}^l)$  \textbf{then}\\
\> \> $s^* := d$; \, $\lambda^*(d, i^*):= \lambda^{k-1}_{d}(i^*)$; \,
$A_d := A_d \cup \{i^*\}$; \, $k_d := k_d + 1$ \\
\>  \>   \textbf{if} \= $k_d < n$  \textbf{then}  
$\begin{bmatrix} \mathbf{w}^{k}_d & \mathbf{r}^{k}_d\end{bmatrix} := 
    \begin{bmatrix} \mathbf{w}^{k-1}_d & \mathbf{r}^{k-1}_d\end{bmatrix}$  \\
\> \textbf{else} \\
\>  \> 
 $(s^*,  i^*) := (s_{d-1}^l, i_{d-1}^l)$ \\
\> \> \textbf{if}  \= $s^{*} = d-1$ \textbf{then} \\
\> \>  \>  $\begin{bmatrix} \mathbf{w}^{k}_d &  \mathbf{r}^{k}_d \end{bmatrix}
:= \begin{bmatrix} \mathbf{w}^{k-1}_d &  \mathbf{r}^{k-1}_d\end{bmatrix} +
     \mathbf{B}({\cdot, i^*}) \begin{bmatrix} 
w^{l-1}_{d-1}(i^*) & r^{l-1}_{d-1}(i^*)\end{bmatrix}$; $A_{d-1} := A_{d-1} \cup \{i^*\}$ \\
\> \> \textbf{else} \{ $s^{*} \leq d-2$ \}\\
\> \>  \> $\begin{bmatrix} w^{k}_{d}(i) & r^{k}_d(i) \end{bmatrix} :=  \begin{bmatrix} 1 &
  R(i) \end{bmatrix} + 
    \sum_{j \in A_{d-1}} b({i, j}) \begin{bmatrix} 
w^{l}_{d-1}(j) & r^{l}_{d-1}(j)\end{bmatrix}, \, i \in \mathbb{X}$  \\
\> \> \textbf{end} \{ if \} \\
\> \> $l := l + 1$ \\
\> \textbf{end} \{ if \} \\
\> $(s_d^k,  j_d^k) := (s^*, i^*)$; \, $k := k + 1$ \\
\textbf{until}  $k_d = n$ \{ repeat \} \\
$L_d := k-1$ 
\end{tabbing}
\end{minipage}}
\end{center}
\label{tab:recalg}
\end{table}

Consider the  algorithm's $d$th stage, $\mathrm{RAG}_{d}$, for a given remaining time $2 \leq d \leq T$, which  is shown in
Table \ref{tab:recalg}, where  $\mathbf{B} = \big(b(i, j)\big)_{i, j \in \mathbb{X}} \triangleq \beta \mathbf{P}$.
The \emph{input} to  $\mathrm{RAG}_{d}$ consists of  
(i) the index values 
$\lambda^*(s, i)$
for smaller horizons $1 \leq s \leq d-1$;  and (ii) sequences  labeled by $l = 1, \ldots, L_{d-1}$ (with $L_{d-1} < (d-1)n$ being part of the input) of augmented states $(s_{d-1}^l, i_{d-1}^l)$, 
and of modified work and reward measures $\mathbf{w}^{l-1}_{d-1} =  \big(w^{l-1}_{d-1}(i)\big)_{i \in \mathbb{X}}$ and $\mathbf{r}^{l-1}_{d-1} =  \big(r^{l-1}_{d-1}(i)\big)_{i \in \mathbb{X}}$
 from the previous  stage.
The \emph{output} of $\mathrm{RAG}_{d}$ gives the
input to the next stage.

Algorithm $\mathrm{RAG}_{d}$ performs $L_d$ steps, labeled by $k = 1,
\ldots, L_d$, to build up the first $L_d$ continuation sets 
$A^{k-1} = (A_1^{k-1}, \ldots, A_d^{k-1})$ of those constructed by algorithm
$\mathrm{AG}({\mathcal{F}_{d}})$ in \ Table \ref{tab:expalg},
using  its input to avoid redundant
computations. 
Note that $|A^{k-1}| = \sum_{s=1}^d |A_s^{k-1}|$, $|A^{k-1}| = k - 1$, and $|A_s^{k-1}| = k_s$ for $1 \leq s \leq d$.
As before,  in the algorithm's notation
$w^{k-1}_s(i)$ and $r^{k-1}_s(i)$ stand for 
$w^{A^{k-1}}_s(i)$ and $r^{A^{k-1}}_s(i)$, respectively.
For a continuation set $A^{k-1}$ as above, we write 
$\widehat{A}^{l-1} = (A_1^{k-1},  \ldots
A_{d-1}^{k-1})$, with $l = k - k_d$. Note that $|\widehat{A}^{l-1}| = l-1$.

The key insight on which the design of algorithm $\mathrm{RAG}_{d}$
is based is that the successive continuation sets
$\widehat{A}^{l-1}$, for $l = 1, 2, \ldots$, corresponding to the $A^{k-1}$ constructed
 by algorithm $\mathrm{RAG}_{d}$, are precisely
those constructed by the previous
stage's algorithm,
$\mathrm{RAG}_{d-1}$.
Exploiting such an insight allows us to 
simplify algorithm $\mathrm{AG}({\mathcal{F}_{d}})$ in Table
\ref{tab:expalg}, as follows.
Consider  step $k$ of algorithm $\mathrm{RAG}_{d}$, which
corresponds to step $k$ of $\mathrm{AG}({\mathcal{F}_{d}})$.
Such a step identifies the augmented state $(s^*,i^*) \in \argmax \big\{\lambda^{k-1}_s{(i)}\colon 1 \leq s \leq d, A_s^{k-1} \subset A_{s+1}^{k-1}, i \in
A_{s+1}^{k-1} \setminus A_s^{k-1}\big\}$ that will be added to the current continuation set $A^{k-1} = (A_1^{k-1}, 
\ldots A_d^{k-1})$ to obtain the next
one, $A^k = A^{k-1} \cup \{(s^*, i^*)\}$, with
the index of augmented state $(s^*, i^*)$ being then given by 
$\lambda^*(s^*, i^*)= \lambda_{s^*}^{{k-1}}(i^*) = r_{s^*}^{k-1}(i^*)/w_{s^*}^{k-1}(i^*)$.
Such a maximization of $\lambda_{s}^{{k-1}}(i)$ is now broken down into two parts: that
for horizon $s = d$, and that for smaller horizons $s < d$, with the first being
the only one that requires actual computations, since the second was already evaluated in previous stages, having as maximizing argument $(s_{d-1}^l, i_{d-1}^l)$, with $l = k - k_d$.
Algorithm $\mathrm{RAG}_{d}$ stops as soon as $A_d^k = \mathbb{X}$, i.e., $k_d = n$, performing $L_d < d n$ steps. 

As the step counter $k$ advances, algorithm $\mathrm{RAG}_{d}$ 
constructs  the
required modified work and reward measures $w_{d}^{k}(i)$
and $r_d^{k}(i)$ using
the recursions obtained in Section 
\ref{s:rmrm}.
The update formulae to use depend on the value of $s^*$. Three cases need to
be considered.
In the first case, $s^* = d$, we use Lemmas \ref{lma:mrmr}(b)
and \ref{lma:mwmr}(b) to conclude that
\[
\begin{bmatrix} r_d^{k}(i) & w_d^{k}(i) \end{bmatrix} =  \begin{bmatrix} r_d^{k-1}(i) &  w_d^{k-1}(i) \end{bmatrix}, \quad i \in \mathbb{X}.
\]
In the second case, $s^* = d-1$, Lemmas \ref{lma:mrmr}(c)
and \ref{lma:mwmr}(c) 
yield that, with  $l = k - k_d$:
\begin{equation}
\label{eq:wrsecc}
\begin{bmatrix} w^{k}_d(i) &  r^{k}_d(i) \end{bmatrix}
= \begin{bmatrix} w^{k-1}_d(i) & r^{k-1}_d(i)\end{bmatrix} +
   b(i, i^*) \begin{bmatrix} 
w^{l-1}_{d-1}(i) &r^{l-1}_{d-1}(i)\end{bmatrix}, \quad i \in
\mathbb{X}.
\end{equation}
Finally, in the third case, $s^* < d-1$, we use Lemmas
\ref{lma:mrmr}(d)
and \ref{lma:mwmr}(d) 
to obtain 
\begin{equation}
\label{eq:wrthirdc}
 \begin{bmatrix} w^{k}_d(i) & r^{k}_d(i) \end{bmatrix} = \begin{bmatrix} 1 &
  R(i) \end{bmatrix} + 
   \sum_{j \in A_{d-1}} b({i,  j}) \begin{bmatrix} 
w^{l}_{d-1}(j) & r^{l}_{d-1}(j)\end{bmatrix}, \quad i \in \mathbb{X},
\end{equation}
where again $l = k - k_d$.
Such recursions are initialized by setting $w^{0}_d(i) \equiv 1$
and $r^{0}_{d}(i) \equiv R(i)$.

The above applies to a stage $2 \leq d \leq T$.
As for the  initial stage $d = 1$, the required quantities are
obtained by setting
$\lambda^*(1, i) = R(i)$,  $w^{l}_1(i) \equiv 1$, and $r^{l}_1(i) \equiv R(i)$.

The following  result  establishes the validity
of the resulting $T$-stage index algorithm RAG, which successively
runs the stages
 $\text{RAG}_1$, $\text{RAG}_2$, \ldots, $\text{RAG}_T$, 
and assesses its time and memory
complexity.
\begin{theorem}
\label{the:mrp}
Algorithm RAG computes all index values
$\big\{\lambda^*(d, i)\colon 1 \leq d \leq T, i \in \mathbb{X}\big\}$  in
$O(T^2 n^3)$ time and $O(T n^2)$
 space.
\end{theorem}
\begin{proof}
The result that algorithm RAG computes index $\lambda^*(d, i)$ is
already proven by the above discussion.
Regarding the time complexity, let us focus on a given
stage $d$ (with $2 \leq d \leq T$), i.e., on algorithm $\text{RAG}_d$.
The description of $\text{RAG}_d$ as given in Table \ref{tab:recalg} shows
that the bulk of the computational work corresponds to the update
in (\ref{eq:wrthirdc}), which entails
\begin{align*}
4 |A_{d-1}| n  +  O(n)  \leq 4 n^2 + O(n)
\end{align*}
AOs. 
Adding up such an upper bound over steps $k
= 1, \ldots, d n$ gives $4 d n^3 + O(d n^2)$ for stage $d$. Then, adding
up over stages $d = 2, \ldots, T$ gives the upper bound
$O(T^2 n^3)$. 

As for memory, the bulk of the storage requirements at
the last, more expensive stage $T$, 
corresponds to 
quantities $\{w^{k}_T(i), r^{k}_T(i)\colon i \in
  \mathbb{X}, 1 \leq k \leq L_T\}$, and $\{w^{l}_{T-1}(i), r^{l}_{T-1}(i)\colon i \in
  \mathbb{X}, 1 \leq l \leq L_{T-1}\}$, with $L_d < d n$, which use $O(T
n^2)$ floating-point locations. Since memory can be reused
from one stage to the next, this gives the overall memory complexity.
\end{proof}

Since  $T n$ index values are  computed, Theorem \ref{the:mrp}
ensures that the \emph{average complexity per index value} of algorithm RAG is
$O(T n^2)$ time and $O(n)$ space.

\subsection{Limited  transitions per state: $O(T^2 n^2)$   index algorithm}
\label{s:iasm}
The $O(T^2 n^3)$ time complexity of the 
 RAG index algorithm holds for a
general project. 
Yet, in many models arising in applications, the  state transition
probability matrix is \emph{sparse}, as only a limited number $N$ of  states, 
which remains
fixed as the total number $n$ of states varies,  can
be reached from any given state. Namely, the following condition
holds. 
\begin{assumption}
\label{ass:ltps}
For every state $i \in \mathbb{X}$, 
$\big|\{j \in \mathbb{X}\colon p({i, j}) > 0\}\big| \leq N$.
\end{assumption}

In such cases, the time complexity of the RAG algorithm
given in the previous section is reduced by an order of magnitude
in the total number $n$ of states.

\begin{proposition}
\label{pro:mrpsl}
Under Assumption $\ref{ass:ltps}$, the 
RAG algorithm  computes all index values
$\big\{\lambda^*(d, i)\colon 1 \leq d \leq T, i \in \mathbb{X}\big\}$
in $O(T^2 n^2)$ time.
\end{proposition}
\begin{proof}
As in Theorem \ref{the:mrp}, 
the bottleneck computation for any  stage $d$ and step $k$
corresponds to the update
in (\ref{eq:wrthirdc}), which now entails no more than
$2 (2 N +1) n  =  O(n)$
AOs. 
Adding up such an upper bound over steps $k
= 1, \ldots, d n$ gives $O(d n^2)$ AOs for stage $d$. Then, adding
up over stages $d = 2, \ldots, T$ gives 
$O(T^2 n^2)$ AOs. 
\end{proof}

\section{Computing the relevant index values of a 
countable-state project}
\label{s:icism}
For a finite-state project,
the RAG algorithm computes index values 
$\lambda^*(d, i)$ for every intermediate horizon $1 \leq d \leq T$ and state $i \in \mathbb{X}$
combination.
Yet,  to deploy the resulting priority-index policy in a
particular instance of the FHMABP
(\ref{eq:dsp}),  the index of each project
need not be evaluated at each such $(d, i)$ pair, but only at the
typically smaller subset of \emph{relevant}
 $(d, i)$, i.e., those that can be reached from the
initial state within the allotted time.
Even for a project  with a countably infinite state space,
provided it satisfies 
Assumption \ref{ass:ltps} --- such as the classic 
Bernoulli bandit  with Beta priors, for which $N = 2$ --- the set of relevant
$(d, i)$ pairs is finite.
 This section presents a modified version of
the RAG algorithm that computes the index values only at
 such relevant $(d, i)$ pairs. 

Consider a project whose state $X(t)$ moves through 
the countable (finite or infinite)
state space $\mathbb{X}$, with its transition probabilities satisfying Assumption \ref{ass:ltps}.
Suppose the project starts at $X(0) = i_0$ with
horizon $T$.
Then, the project's index $\lambda^*(d, i)$ need be
evaluated only at pairs $(d, i)$ such that: 
(i) $1 \leq d \leq T$; and (ii) 
  state $i$ belongs to the finite set $\mathbb{X}_{T-d}(i_0)$ of states that 
 can be reached within $T-d$  periods starting at
 $i_0$. Note that $n_d(i_0) \triangleq |\mathbb{X}_{T-d}(i_0)|  \leq \sum_{s=0}^{T-d} N^s$, and that
\[
\{i_0\} = \mathbb{X}_0(i_0) \subseteq \mathbb{X}_1(i_0) \subseteq
\cdots \subseteq \mathbb{X}_T(i_0).
\]
Let $\mathbb{Y}_T^{\{0, 1\}}(i_0) \triangleq \big\{(d, i)\colon 1 \leq d
\leq T, i \in \mathbb{X}_{T-d}(i_0)\big\}$ be the finite set of such
 relevant $(d, i)$ pairs, and
let $\mathbb{Y}_T(i_0) \triangleq \mathbb{Y}_T^{\{0, 1\}}(i_0)
\cup \{(0, \Omega)\}$. Note that $\mathbb{Y}_T^{\{0, 1\}}(i_0)$ contains 
$L_T(i_0) \triangleq |\mathbb{Y}_T^{\{0, 1\}}(i_0)| = \sum_{d=1}^T n_d(i_0)$ augmented states $(d, i)$.

Consider now  an
 auxiliary finite-state infinite-horizon project, whose state 
$Y(t)$ evolves over time periods $t \geq 0$ through the state space 
$\mathbb{Y}_T(i_0)$, defined as in Section \ref{s:icpclagia} but using
$\mathbb{Y}_T(i_0)$ and $\mathbb{Y}_T^{\{0, 1\}}(i_0)$ in place of
$\mathbb{Y}_T$ and $\mathbb{Y}_T^{\{0, 1\}}$.
Now,  let $G'(d, i)$ be the 
\emph{modified Gittins index} of such an auxiliary project, as defined by (\ref{eq:gsi}). 
The interest of introducing such an auxiliary project and its modified
Gittins index $G'(d, i)$ is, as in Section \ref{s:icpclagia}, 
that the latter is precisely the
finite-horizon index $\lambda^*(d, i)$.
The proof of the next result follows along the same lines as that
of Proposition \ref{pro:rgi}, and is hence omitted.

\begin{proposition}
\label{pro:rgi2}
$\lambda^*(d, i) = G'(d, i)$ for $(d, i) \in  \mathbb{Y}_T^{\{0, 1\}}(i_0)$.
\end{proposition}

As in Section \ref{s:agia}, 
Proposition \ref{pro:rgi2} allows us to obtain the finite set of relevant index 
values $\lambda^*(d, i)$ for $(d, i) \in  \mathbb{Y}_T^{\{0, 1\}}(i_0)$
by running the adaptive-greedy 
algorithm $\mathrm{AG}({\mathcal{F}_T})$ in Table \ref{tab:agaf}.
Yet, note that the definition of active-set family 
$\mathcal{F}_T$ in (\ref{eq:fhat}) used in 
Section  \ref{s:agia} 
must be modified in the present setting to 
$\mathcal{F}_T \triangleq 2^{\mathbb{Y}_T^{\{0, 1\}}(i_0)}$. 

We can now use the 
recursionspresented in Section
\ref{s:rmrm},  along the lines in Section \ref{s:ria}, to reformulate the one-pass algorithm $\mathrm{AG}({\mathcal{F}_T})$ into a 
$T$-stage recursive algorithm, which we denote by 
$\mathrm{RAG}(i_0)$ to emphasize its dependence on the initial project 
state $i_0$.
Table \ref{tab:recalginf} shows stage $d$ of such an algorithm, 
which we denote by $\mathrm{RAG}_d(i_0)$, 
for $2 \leq d \leq T$. 

\begin{table}[htb]
\caption{Stage $2 \leq d \leq T$ of  
index algorithm $\mathrm{RAG}(i_0)$.}
\begin{center}
\fbox{%
\begin{minipage}{6in}
\textbf{ALGORITHM} $\mathrm{RAG}_{d}(i_0)$ \\
\textbf{Input:}
$\{\lambda^*(s, i)\colon 1 \leq s \leq d-1, i \in \mathbb{X}_{T-s}(i_0)\}, \{(s_{d-1}^l, i_{d-1}^l), \mathbf{w}^{l-1}_{d-1}, \mathbf{r}^{l-1}_{d-1}\colon 1 \leq l \leq
  L_{d-1}\}$
 \\
\textbf{Output:}
$\{\lambda^*(s, i)\colon 1 \leq s \leq d, i \in \mathbb{X}_{T-s}(i_0)\}, 
 \{(s_{d}^k, i_{d}^k), \mathbf{w}^{k-1}_d, \mathbf{r}^{k-1}_d\colon 1 \leq k \leq L_d\}$ 
\begin{tabbing}
$k := 1$; \, $l := 1$; \, $A_{d-1} := A_d := \emptyset$; \,   
$k_d := 0$; \, 
$\begin{bmatrix} \mathbf{w}^{0}_{d} & \mathbf{r}^{0}_d\end{bmatrix} :=  \begin{bmatrix} \mathbf{1} &
  \mathbf{R} \end{bmatrix}$ \\
\textbf{repeat} \=    \{ note: $A^{k-1} = (A_1^{k-1}, \ldots, A_d^{k-1}) = (A_1, \ldots, A_d)$ \} \\
\> $\lambda^{k-1}_{d}(i) := r^{k-1}_{d}(i)/w^{k-1}_{d}(i),  i \in \mathbb{X}_{T-d}(i_0) \setminus A_d$ \\
\>  \textbf{pick} 
 $\displaystyle i^{*} \in \argmax
\big\{\lambda^{k-1}_{d}(i)\colon i \in \mathbb{X}_{T-d}(i_0) \setminus A_d\big\}$ \\
\> \textbf{if} \= $\lambda^{k-1}_{d}(i^*) \geq \lambda^*(s_{d-1}^l, i_{d-1}^l)$  \textbf{then}\\
\> \> $s^* := d$; \, $\lambda^*(d, i^*):= \lambda^{k-1}_{d}(i^*)$; \,
$A_d := A_d \cup \{i^*\}$; \, $k_d := k_d + 1$ \\
\>  \>   \textbf{if} \= $k_d < n_d(i_0)$  \textbf{then}  
$\begin{bmatrix} w^{k}_d(i) & r^{k}_d(i)\end{bmatrix} := 
    \begin{bmatrix} w^{k-1}_d(i) & r^{k-1}_d(i)\end{bmatrix}, \, i \in \mathbb{X}_{T-d}(i_0)$  \\
\> \textbf{else} \\
\>  \> 
 $(s^*,  i^*) := (s_{d-1}^l, i_{d-1}^l)$ \\
\> \> \textbf{if}  \= $s^{*} = d-1$ \textbf{then} \\
\> \>  \>  $\begin{bmatrix} w^{k}_d(i) &  r^{k}_d(i) \end{bmatrix}
:= \begin{bmatrix} w^{k-1}_d(i) &  r^{k-1}_d(i)\end{bmatrix} +
     b({i, i^*}) \begin{bmatrix} 
w^{l-1}_{d-1}(i^*) & r^{l-1}_{d-1}(i^*)\end{bmatrix}, i \in \mathbb{X}_{T-d}(i_0)$ \\
\> \>  \>   $A_{d-1} := A_{d-1} \cup \{i^*\}$ \\
\> \> \textbf{else} \{ $s^{*} \leq d-2$ \}\\
\> \>  \> $\begin{bmatrix} w^{k}_{d}(i) & r^{k}_d(i) \end{bmatrix} :=  \begin{bmatrix} 1 &
  R(i) \end{bmatrix} + 
    \sum_{j \in A_{d-1}} b({i, j}) \begin{bmatrix} 
w^{l}_{d-1}(j) & r^{l}_{d-1}(j)\end{bmatrix}, \, i \in \mathbb{X}_{T-d}(i_0)$  \\
\> \> \textbf{end} \{ if \} \\
\> \> $l := l + 1$ \\
\> \textbf{end} \{ if \} \\
\> $(s_d^k,  j_d^k) := (s^*, i^*)$; \, $k := k + 1$ \\
\textbf{until}  $k_d = n_d(i_0)$ \{ repeat \} \\
$L_d := k-1$ 
\end{tabbing}
\end{minipage}}
\end{center}
\label{tab:recalginf}
\end{table}

The next  result  ensures the validity
of the resulting $T$-stage index algorithm $\mathrm{RAG}(i_0)$, 
which successively runs 
the stages $\text{RAG}_1(i_0)$, $\text{RAG}_2(i_0)$, \ldots,
$\text{RAG}_T(i_0)$, and assesses its complexity.
For the latter purpose, we further assume a quadratic growth rate for 
$|\mathbb{X}_s(i_0)|$ in the remaining time $s$, as this is common in applications.

\begin{assumption}
\label{ass:xti0growth} 
$|\mathbb{X}_s(i_0)| = O(s^2)$.
\end{assumption}

Under Assumption  \ref{ass:xti0growth}, algorithm $\mathrm{RAG}(i_0)$ computes  
$L_T(i_0) = O(T^3)$ index values $\lambda^*(d, i)$.

\begin{theorem}
\label{the:mrpsp}
Algorithm $\mathrm{RAG}(i_0)$ computes index values
$\big\{\lambda^*(d, i)\colon (d, i) \in \mathbb{Y}_T^{\{0,
  1\}}(i_0)\big\}$ in $O(T^6)$ time
and $O(T^5)$ space.
\end{theorem}
\begin{proof}
The result that algorithm $\mathrm{RAG}(i_0)$ computes the stated
values of index 
$\lambda^*(d, i)$ is
already proven by the above discussion.
Regarding the time complexity, let us focus on a given
stage $d$ (with $2 \leq d \leq T$), i.e., on algorithm $\text{RAG}_d(i_0)$.
The description of $\text{RAG}_d(i_0)$ as given in Table
\ref{tab:recalginf}
 shows
that the computational bottleneck is the update
\[
\begin{bmatrix} w^{k}_{d}(i) & r^{k}_d(i) \end{bmatrix} :=  \begin{bmatrix} 1 &
  R(i) \end{bmatrix} + 
    \sum_{j \in A_{d-1}} b({i, j}) \begin{bmatrix} 
w^{l}_{d-1}(j) & r^{l}_{d-1}(j)\end{bmatrix}, \, i \in \mathbb{X}_{T-d}(i_0),
\] 
which entails no more than
$2 (2 N + 1) n_d(i_0)$
AOs.
Adding up such an upper bound over steps $k
= 1, \ldots, L_d \leq L_d(i_0)$ gives no more than $2 (2 N + 1) n_d(i_0) L_d(i_0)$ for stage $d$. Then, adding
up over stages $d = 2, \ldots, T$ and using Assumption 
\ref{ass:xti0growth}, gives an upper bound of $O(T^6)$ AOs.

As for memory, the bulk of storage at
 stage $d$ 
corresponds to 
the $w^{k-1}_{s}(i)$ and  $r^{k-1}_{s}(i)$ for $1 \leq k \leq L_s(i_0)$,
   $i \in \mathbb{X}_{T-s}(i_0)$, and $s = d-1, d$. Now, for every $s$, 
$2 n_s(i_0) L_s(i_0) = O(T^5)$ (see Assumption
\ref{ass:xti0growth}) floating-point memory locations are needed. 
\end{proof}

To illustrate, in the case of  the Bernoulli bandit
model with Beta priors, where the  state is a pair $(i, j) \in \{1, 2,
\ldots\}^2$ giving the parameters of the corresponding posterior
Beta distribution, suppose one wants to compute the index 
values $\lambda^*\big(d, (i, j)\big)$ for states $(i, j)$ that can be reached from
a given initial state $(i_0, j_0)$ within $T$ periods. 
For such a model, Assumption \ref{ass:ltps} holds with $N = 2$ and, since
\[
\mathbb{X}_s(i_0, j_0) =
\big\{(i, j) \geq (i_0, j_0)\colon (i-i_0) + (j-j_0) \leq s\big\},
\]
Assumption \ref{ass:xti0growth} also holds, since 
$|\mathbb{X}_s(i_0, j_0)| = 1 + \cdots + (s+1) = (s+1) (s+2) / 2 =
O(s^2)$.
Hence,  the total number of index 
values that is computed by algorithm $\text{RAG}(i_0, j_0)$ is 
$L_T(i_0, j_0) = T (T+1) (T+2) / 6 = O(T^3)$.

The reader may wonder how the complexity 
results in Theorem \ref{the:mrpsp}, as applied to the Bernoulli bandit
model with Beta priors, 
compare with those reported
in \citet[p.\ 139]{gi89}, which might appear better at first glance, 
being $O(T^4)$ AOs and $O(T^2)$ memory.
 The answer is that both complexity counts cannot be
 meaningfully compared, for the following reasons: (i) the purpose of the algorithm in \citet[Sec.\ 7]{gi79} is
to
approximate the Gittins index $\lambda^*(i_0, j_0)$ at a single state $(i_0, j_0)$ by
$\lambda^*\big(T, (i_0, j_0)\big)$, for which only the subset of $1 + \cdots + T =
O(T^2)$ index values of the form $\lambda^*\big(d, (i, j)\big)$, for $1 \leq d
\leq T$ and $(i-i_0) + (j-j_0) = T-d$, needs to be evaluated; and (ii) 
the Gittins algorithm calls for solving a continuous optimization problem
of the form (\ref{eq:lambdadtaulamb}) at each step, which is not an 
elementary operation, whereas the RAG algorithm herein
performs only  arithmetic operations.

\section{Block implementation of the  RAG index algorithm}
\label{s:riraga}
This section presents an efficient implementation of the RAG index algorithm.
A naive implementation
 which directly codes the algorithm's
update formulae will be found to be rather slow and inefficient, even 
for instances with a moderate number of states or horizon.
The reason is that the bottleneck computation, which is the update in (\ref{eq:wrthirdc}),
 involves repeated 
multiplications of large matrices with \emph{noncontiguous
  memory-access patterns}, requiring expensive 
 \emph{gather} and \emph{scatter} memory operations.
Such patterns cause severe 
inefficiencies in linear algebra algorithms,
 due to the mismatch between the speeds of
processors (fast) and of memory access (slow) in contemporary
computers.
The main approach to 
 reduce such inefficiencies, exploiting both vectorization and parallelism features of advanced computer architectures, is to design \emph{block implementations}. These 
aim to maximize the arithmetic operations performed per memory
access, by rearranging bottleneck computations as linear algebra operations on \emph{contiguous blocks} of data (e.g., matrix-matrix multiplications), thus attaining a sort of economies of scale in computation. 
See  \cite{dongeijk00}.

\begin{table}[htb]
\caption{Block implementation of stage $2 \leq d \leq T$ of the RAG index algorithm.}
\begin{center}
\fbox{%
\begin{minipage}{6in}
\textbf{ALGORITHM} $\mathrm{BlockRAG}_{d}$ \\
\textbf{Input:}
$\big\{\lambda^*(s, i)\big\}_{1 \leq s \leq d-1, i \in \mathbb{X}}, \big\{(s_{d-1}^l, i_{d-1}^l), \widehat{\mathbf{w}}^{l-1}_{d-1}, \widehat{\mathbf{r}}^{l-1}_{d-1}\}_{1 \leq l \leq L_{d-1}}$, $\mathbf{w}_{d-1}^*, \mathbf{r}_{d-1}^*$
 \\
\textbf{Output:}
$\big\{\lambda^*(s, i)\big\}_{1 \leq s \leq d, i \in \mathbb{X}}, \big\{(s_{d}^k, i_{d}^k), \widehat{\mathbf{w}}^{k-1}_d, \widehat{\mathbf{r}}^{k-1}_d\big\}_{1 \leq k \leq L_d}$, $\mathbf{w}_{d}^*, \mathbf{r}_{d}^*$
\begin{tabbing}
$k := 1$; \, $l := 1$; \, $A_d := \emptyset$; \,   
$k_d := 0$; \, 
$\begin{bmatrix} \widehat{\mathbf{w}}^{0}_{d} & \widehat{\mathbf{r}}^{0}_d\end{bmatrix} :=  \begin{bmatrix} \mathbf{1} &
  \mathbf{R} \end{bmatrix}$ \\
\textbf{repeat} \=   \{ note: $A^{k-1} = (A_1^{k-1}, \ldots, A_d^{k-1}) = (A_1, \ldots, A_d)$ \} \\
\> $\lambda^{k-1}_{d}(i) := r^{k-1}_{d}(i)/w^{k-1}_{d}(i),  i \in \mathbb{X} \setminus A_d^{k-1}$;  \textbf{pick} 
 $\displaystyle i^{*} \in \argmax
\big\{\lambda^{k-1}_{d}(i)\colon i \in \mathbb{X} \setminus A_d\big\}$ \\
\> \textbf{if} \= $\lambda^{k-1}_{d}(i^*) \geq \lambda^*(s_{d-1}^l, i_{d-1}^l)$  \textbf{then}\\
\> \> $s^* := d$;  $\lambda^*(d, i^*):= \lambda^{k-1}_{d}(i^*)$; 
$A_d := A_d \cup \{i^*\}$;  $k_d := k_d + 1$;  $k^*(i^*) := k-1$ \\
\>  \>   \textbf{if} \= $k_d < n$  \textbf{then}  \\
\>  \>   \> 
$\begin{bmatrix}w_d(i^*) & r_d(i^*) \end{bmatrix} :=  \begin{bmatrix}w^{k-1}_d(i^*) &  r^{k-1}_d(i^*) \end{bmatrix}$; \,
$\begin{bmatrix} \widehat{\mathbf{w}}^{k}_d & \widehat{\mathbf{r}}^{k}_d\end{bmatrix} := 
    \begin{bmatrix} \widehat{\mathbf{w}}^{k-1}_d & \widehat{\mathbf{r}}^{k-1}_d\end{bmatrix}$  \\
\>  \>    \textbf{end } \{ if \} \\
\> \textbf{else} \\
\>  \> 
 $(s^*,  i^*) := (s_{d-1}^l, i_{d-1}^l)$ \\
\> \> \textbf{if}  \= $s^{*} = d-1$ \textbf{then} \\
\> \>  \>  $\begin{bmatrix} \widehat{\mathbf{w}}^{k}_d &  \widehat{\mathbf{r}}^{k}_d \end{bmatrix}
:= \begin{bmatrix} \widehat{\mathbf{w}}^{k-1}_d &  \widehat{\mathbf{r}}^{k-1}_d\end{bmatrix} +
     \mathbf{B}({\cdot, i^*}) \begin{bmatrix} 
w_{d-1}(i^*) & r_{d-1}(i^*)\end{bmatrix}$ \\
\> \> \textbf{else} \{ $s^{*} \leq d-2$ \}\\
\> \> \> $\begin{bmatrix} \widehat{\mathbf{w}}^{k}_{d} & \widehat{\mathbf{r}}^{k}_d\end{bmatrix} :=  \begin{bmatrix} \mathbf{1} &
  \mathbf{R} \end{bmatrix} + 
    \begin{bmatrix} 
\widehat{\mathbf{w}}^{l}_{d-1} & \widehat{\mathbf{r}}^{l}_{d-1}\end{bmatrix}$ \\
\> \> \textbf{end} \{ if \} \\
\> \> $l := l + 1$ \\
\> \textbf{end} \{ if \} \\
\> $(s_d^k,  j_d^k) := (s^*, i^*)$; \, $k := k + 1$ \\
\textbf{until}  $k_d = n$ \{ repeat \} \\
$L_d := k-1$ \\
\textbf{for} \= $i \in \mathbb{X}$ \textbf{do}:
$\widehat{w}^{k'}_d(i) := \widehat{r}^{k'}_d(i) := 0$, \, $k' = 0,  \ldots, k^*(i)$; 
\textbf{end} \{ for \} \\
$\begin{bmatrix} \widehat{\mathbf{w}}^{k'-1}_d & \widehat{\mathbf{r}}^{k'-1}_d\end{bmatrix}_{k' = 1}^{L_d} := 
    \mathbf{B} \begin{bmatrix} \widehat{\mathbf{w}}^{k'-1}_d & \widehat{\mathbf{r}}^{k'-1}_d\end{bmatrix}_{k' = 1}^{L_d}$
\end{tabbing}
\end{minipage}}
\end{center}
\label{tab:recalgblock}
\end{table}

Table \ref{tab:recalgblock} presents a block implementation of stage $d$ of the RAG algorithm, denoted by 
$\mathrm{BlockRAG}_{d}$. 
The input to $\mathrm{BlockRAG}_{d}$ differs from that to  $\mathrm{RAG}_{d}$ in that (i) it takes a matrix of vectors
$\widehat{\mathbf{w}}_{d-1}^{l-1}$ and $\widehat{\mathbf{r}}_{d-1}^{l-1}$  instead of the $\mathbf{w}_{d-1}^{l-1}$ and $\mathbf{r}_{d-1}^{l-1}$, where $\widehat{w}_{d-1}^{l-1}(i) \triangleq \sum_{j \in A_{d-1}^{l-1}} b(i, j) w_{d-1}^{l-1}(j)$
and $\widehat{r}_{d-1}^{l-1}(i) \triangleq \sum_{j \in A_{d-1}^{l-1}} b(i, j) r_{d-1}^{l-1}(j)$,  where $A_{d-1}^{l-1}$ corresponds to the 
continuation sets $A^{l-1} = (A_1^{l-1}, \ldots, A_{d-1}^{l-1})$ generated by algorithm $\mathrm{BlockRAG}_{d-1}$; and (ii) it incorporates 
vectors $\mathbf{w}_{d-1}^* = \big(w_{d-1}^*(i)\big)_{i \in \mathbb{X}}$ and $\mathbf{r}_{d-1}^* = \big(r_{d-1}^*(i)\big)_{i \in \mathbb{X}}$, 
where $w_{d-1}^*(i) \triangleq w_{d-1}^{l-1}(i)$ and $r_{d-1}^*(i) \triangleq r_{d-1}^{l-1}(i)$, with $A_{d-1}^l$ being the first continuation set for horizon $d-1$ in 
algorithm $\mathrm{BlockRAG}_{d-1}$ that
contains state $i$.
The output of algorithm $\mathrm{BlockRAG}_{d}$ differs accordingly from that of $\mathrm{RAG}_{d}$.

Algorithm $\mathrm{BlockRAG}_{d}$ implements the stage-$d$ update in (\ref{eq:wrthirdc}) as  $\begin{bmatrix} \widehat{\mathbf{w}}^{k}_{d} & \widehat{\mathbf{r}}^{k}_d\end{bmatrix} :=  \begin{bmatrix} \mathbf{1} &
  \mathbf{R} \end{bmatrix} + 
    \begin{bmatrix} 
\widehat{\mathbf{w}}^{l}_{d-1} & \widehat{\mathbf{r}}^{l}_{d-1}\end{bmatrix}$ (note that, at this point in the algorithm, $\begin{bmatrix} \widehat{\mathbf{w}}^{k}_{d} & \widehat{\mathbf{r}}^{k}_d\end{bmatrix} = \begin{bmatrix} \mathbf{w}^{k}_{d} & \mathbf{r}^{k}_d\end{bmatrix}$).
As for the bottleneck computation, it has been moved out of the loop, to the last line in Table \ref{tab:recalgblock}, as a block matrix-matrix multiplication that computes the
$\begin{bmatrix} \widehat{\mathbf{w}}^{k}_{d} & \widehat{\mathbf{r}}^{k}_d\end{bmatrix}$ to be used in the next stage.

A straightforward modification of algorithm $\mathrm{BlockRAG}_{d}$ gives a corresponding block implementation of 
algorithm $\mathrm{RAG}_d(i_0)$ (see Table \ref{tab:recalginf}), which we denote by $\mathrm{BlockRAG}_{d}(i_0)$.
The author has coded in Fortran the resulting block algorithms $\mathrm{BlockRAG}$ and $\mathrm{BlockRAG}(i_0)$, which have been used
in the experiments reported in Section \ref{s:ce}. 

\section{Computational experiments}
\label{s:ce}
This section reports the results of a computational study, based on the author's Fortran implementations of the algorithms discussed above (which are available for
download at http://alum.mit.edu/www/jnimora, under the link ``Original software codes"), which benchmarks the actual runtime and memory performance of the proposed RAG algorithm against the
calibration method, and fits the measured performance to the theoretical complexity.
\subsection{Index computation for finite-state projects}
\label{s:rtmua}
This experiment 
benchmarks the RAG index algorithm against the calibration method 
on finite-state projects, measuring actual runtime performance
and storage requirements.
Recall that Theorem \ref{the:mrp}
establishes an
$O(T^2 n^3)$ time complexity and an $O(T^2 n^2)$ memory complexity
for the RAG algorithm
 on an $n$-state $T$-horizon project, whereas Proposition
   \ref{pro:caaoc}(b) shows that the time and memory complexities 
 for the calibration
method, when used with a grid of $L$ $\lambda$-values, are $O(L T n^2)$ and 
$O(L T n)$, respectively.

The experiment uses the block implementations (see 
Sections \ref{s:caaic} and
  \ref{s:riraga}) designed and coded in Fortran by
  the author  of the calibration
method and the 
RAG algorithm.
The codes
 were compiled using the latest release
at the time of writing of the Intel Visual Fortran Compiler
Professional Ed.\ 11.1  (Update 6). 
Such implementations use high-performance threaded
 routines from the Intel Math Kernel Library for
bottleneck computations (in particular the BLAS Level 3 
\verb+DGEMM+ subroutine 
 for matrix-matrix multiplication), which can harness to a substantial extent
 the parallel processing power of the platform employed: an 
HP z800 workstation with two quad-core 3.33 GHz
Intel Xeon processors w5590 and 48 GB of memory,
under Windows 7 x64.
Both methods were tested on 20 project instances, 
with state-space sizes $n = 100,
200, \ldots, 2000$, and a horizon of $50$. 
The transition probability matrix of the $n$-state instance
was obtained by scaling an  $n \times n$ matrix with pseudorandom
$\text{Uniform}(0, 1)$ entries, dividing each row by its sum.
Immediate rewards were also drawn from a pseudorandom $\text{Uniform}(0, 1)$
distribution.
The discount factor used was $\beta = 1$.
For each instance, the index values $\lambda^*(T, i)$ 
for $1 \leq i \leq n$ and $1 \leq T \leq 50$ 
 were evaluated using the RAG algorithm and the
calibration method, the latter for 3, 4, and 5 significant digits of
accuracy (partitioning the unit interval $[0, 1]$
with a grid of $L = 10^m + 1$ equally-spaced $\lambda$-values, for $m = 3, 4, 5$).
For each method and 
intermediate horizon $T = 1, \ldots, 50$,
 the  wall-clock cumulative runtime $y(T, n)$  to compute the index values 
$\lambda^*(d, i)$ for $1 \leq d \leq T$ and $i = 1, \ldots, n$ was
measured using the  Fortran intrinsic subroutine \verb+system_clock+.

\begin{figure}[ht!]
\centering
\includegraphics[height=3in,width=6in,keepaspectratio]{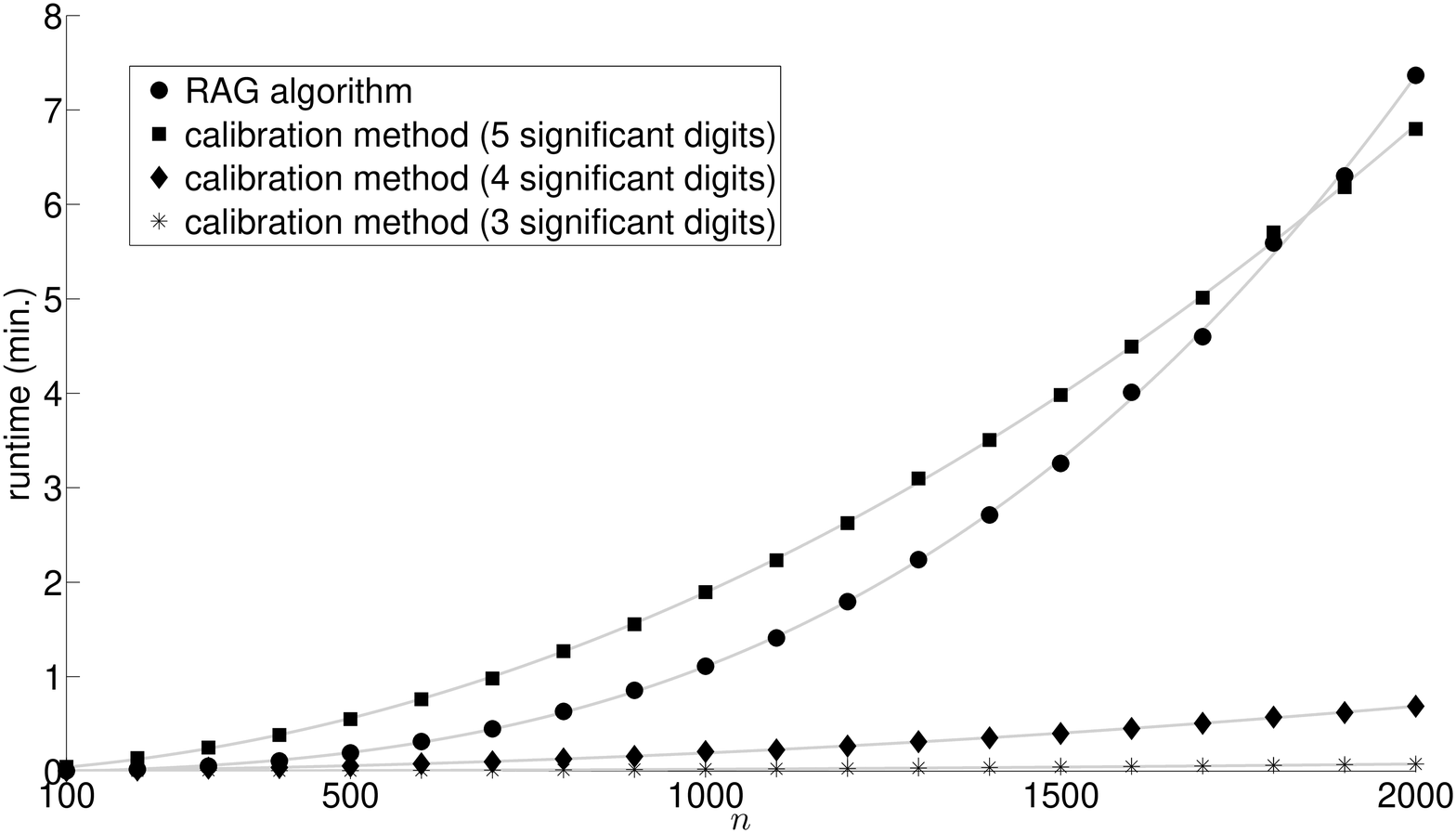}
\caption{{\bf Runtimes vs.\
    number of states $n$ for horizon $T = 50$}.}
\label{plot1Timingexpb}
\end{figure}

Figure \ref{plot1Timingexpb} plots 
the recorded cumulative runtimes (in minutes) versus the number $n$ of states for horizon $T = 50$. 
The gray solid lines shown are 
  polynomial least-squares (LS) fits for the predicted runtimes, of 3rd order for the RAG algorithm, 
and of 2nd order for the calibration method, corresponding to the 
theoretical complexities. 
In the case of the RAG algorithm, the 3rd-order LS fit $\widehat{y}(n, 50)$ for the
predicted runtime  of an instance with $n$ states and 
horizon $50$ is  
 $\widehat{y}(n, 50)=  10^{-10} (7.82 \, n^3 + 2.18 \times 10^3 \, n^2
 +1.26
 \times 10^6 \, n - 1.89 10^8)$. 
To measure of the quality of fit, we use the  \emph{root mean square error}
(RMSE). In this case, the RMSE is $0.04$ min., which indicates that the fit is rather tight, considering the range of runtimes. 
To assess the validity of the theoretical cubic complexity on
$n$, the data were also fitted by polynomials of one order less
and of one order more than 3. The 4th-order polynomial fit has a
spurious negative leading coefficient $-10^{-10} \times 1.67$, with its 
RMSE being about the same as that for the 3rd-order fit. 
As for the 2nd-order LS fit, the RMSE  degrades
significantly, to 
$0.13$ min. These results show that the 3rd-order polynomial gives the best fit.
Note further that, despite its higher complexity, the RAG algorithm
is actually faster than the calibration method with 5 significant
digits up to and including
$n = 1800$ states.

\begin{figure}[ht!]
\centering
\includegraphics[height=3in,width=6in,keepaspectratio]{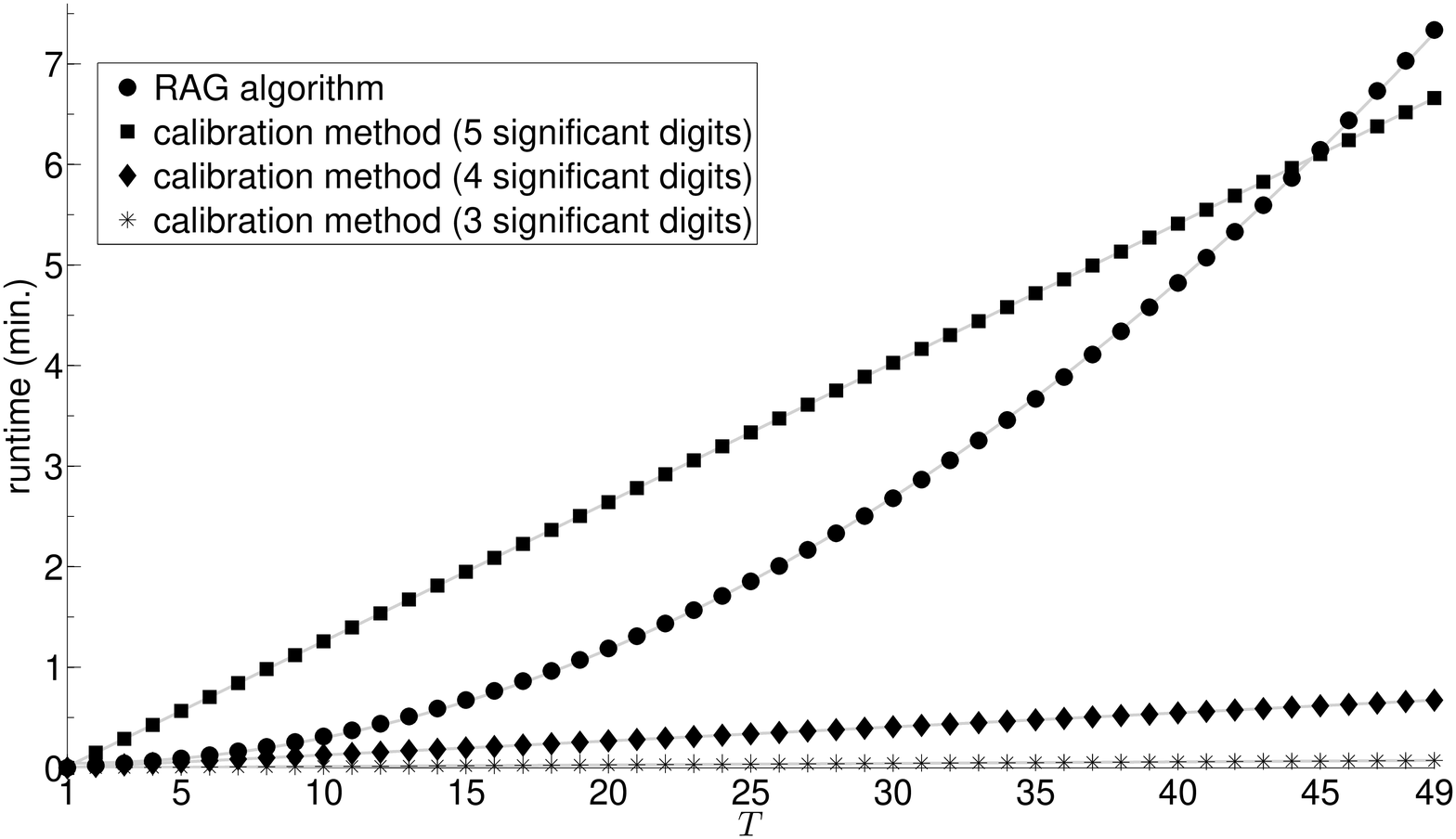}
\caption{{\bf Cumulative runtimes vs.\ horizon (stage) $T$ for $n = 2000$ states.}}
\label{plot2Timingexpb}
\end{figure}

Figure \ref{plot2Timingexpb} plots 
the measured cumulative runtimes versus the intermediate horizon or stage $T = 1, 2, \ldots, 49$
 for the instance with $n = 2000$ states. Cumulative runtimes for the final stage $T = 50$ are not
 included, since the RAG algorithm does not perform the
 bottleneck update at the last stage, and hence the latter's cumulative runtime is about the
 same as that for the previous stage $T = 49$.
The  solid lines shown are 
  polynomial LS fits for the predicted cumulative runtimes, of 2nd order for the RAG algorithm, 
and of 1st order (linear fit) for the calibration method, as predicted by the 
theoretical complexities. 
In the case of the RAG algorithm, the 2nd-order LS fit $\widehat{y}(2000, T)$ for the
predicted cumulative runtime  of a $2000$-state instance up to and including stage $T$ 
  is  
 $\widehat{y}(2000, T)=  10^{-3} (3.16 \, T^2 -6.68 \, T + 49.14)$. The RMSE
 is $0.015$ min., slightly under 1 sec., a very small value relative to the range of runtimes. 
To test the validity of the theoretical quadratic complexity on
$T$, the data were also fitted by a polynomial with one order less
and of one order more than 2. While the linear LS fit is clearly inadequate,
using a polynomial of order 3 gives the 
predicted cumulative runtime fit $\widehat{y}(2000, T)=  10^{-3} (0.0063 \, T^3 + 2.68 \, T^2 + 2.89 \, T + 7.25 )$, with the
RMSE dropping to $0.003$ min. 
Since the RMSE for the 2nd-order fit is already very small, and the leading coefficient of the 3rd-order fit is rather small, we conclude that the 
cumulative runtime performance is best fittedted by a 2nd-order polynomial, consistently with the theoretical complexity in $T$. 
Still, despite its higher complexity, the RAG algorithm
is faster than the calibration method with 5 significant digits up to and including
a horizon of
$T = 44$.

\begin{figure}[ht!]
\centering
\includegraphics[height=3in,width=6in,keepaspectratio]{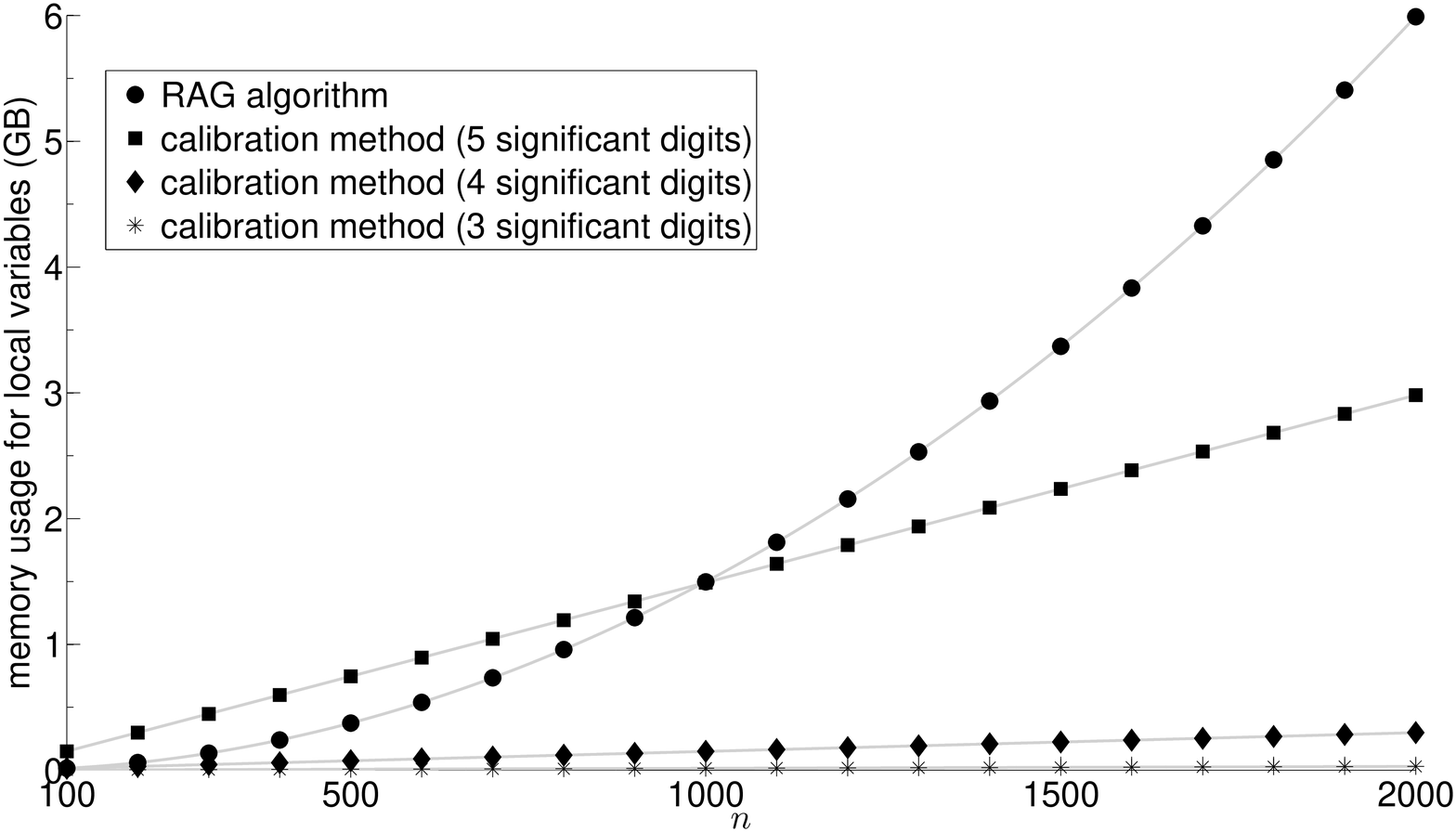}
\caption{{\bf Memory usage vs.\ number of states for horizon $T = 50$.}}
\label{plot3Timingexp}
\end{figure}

Figure \ref{plot3Timingexp} plots 
the required memory storage for floating-point local variables (in GB, excluding input and
  output, and using 8-byte double-precision numbers)  versus $n$ {\bf for $T = 50$}. 
In the author's implementations, the RAG algorithm uses $(4 T + 1) n^2 + 5 n$ floating-point storage locations for local variables, whereas the 
calibration method with a grid of size $L$ uses $2 L n + L$ locations.
Note that, despite its higher complexity, the RAG algorithm
uses less memory than the calibration method with 5 significant digits up to and including 
$n = 900$ states. Note further that the local memory storage of the RAG algorithm grows linearly in the horizon $T$, whereas that of the calibration  method remains constant as
$T$ varies.

\subsection{Index computation for an infinite-state project}
\label{s:rtmua2}
The next experiment aims to assess the actual runtime and memory
 performance of the 
index algorithm in Section \ref{s:icism}
for a countably-infinite state project, for which the classic Bernoulli bandit
model with Beta priors is chosen,  by benchmarking it against
the
calibration method. Recall from Section \ref{s:icism} that Theorem \ref{the:mrpsp}
establishes an
$O(T^6)$ time complexity and an $O(T^5)$ memory complexity
for such a version of the RAG algorithm, which in the case of a
Bernoulli bandit starting at $(i_0, j_0)$ with a horizon $T$ computes $T (T+1) (T+2) / 6$
relevant index values.
For fairness of comparison, the calibration method was modified to compute
approximately only such relevant index values.
To improve runtimes and exploit the reduced arithmetic and memory 
operations due to sparsity of the transition probability matrix, the
author
developed 
Fortran implementations that use  threaded
 routines from the Intel Math Kernel Library, 
in particular the Sparse BLAS Level 3 
\verb+MKL_DCOOMM+ subroutine 
 for sparse matrix multiplication.

Taking $(1, 1)$ as the initial state, the
algorithm $\mathrm{RAG}(1, 1)$ in Section \ref{s:icism} was run on
 instances with horizons $T = 20, 25, \ldots, 90$, 
computing all relevant index values in each case. 
As before, 
the calibration method (with $3$ to $5$ significant digits) was used
for comparison.

\begin{figure}[ht!]
\centering
\includegraphics[height=3in,width=6in,keepaspectratio]{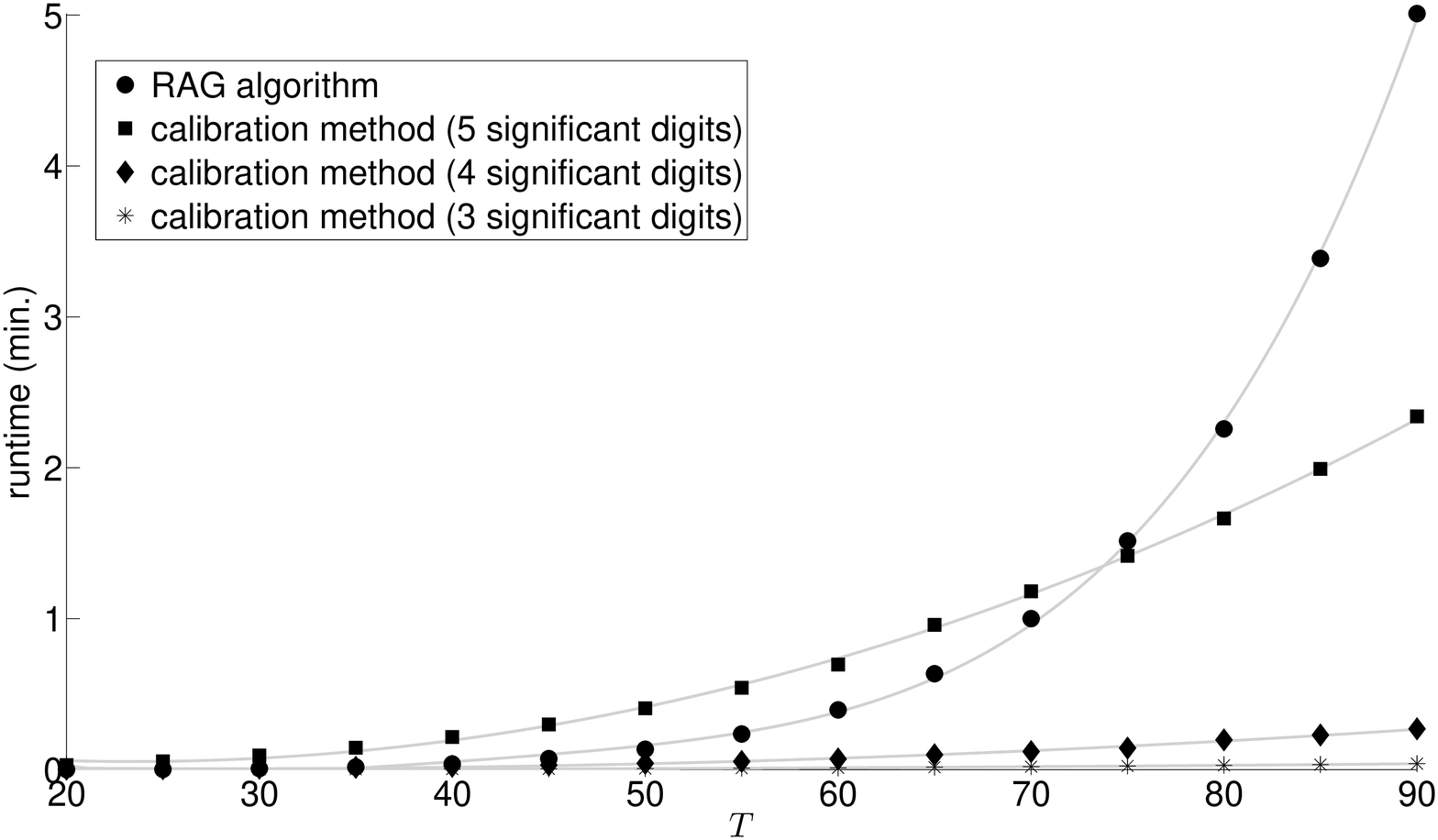}
\caption{Runtimes  vs.\ horizon $T$.}
\label{plot1BernExpb}
\end{figure}

Figure \ref{plot1BernExpb} plots 
the measured runtimes versus the horizon $T$. 
The  solid lines shown were obtained by 
  polynomial least-squares (LS) fit, of order $4$ for the $\mathrm{RAG}(1, 1)$ algorithm, 
and of order $2$ for the calibration method. 
In the case of the $\mathrm{RAG}(1, 1)$ algorithm, the 6th-order LS fit for the
predicted runtime $\widehat{y}(T)$ of an instance with $T$ remaining periods is  
 $\widehat{y}(T)=  10^{-10} (1.71 \, T^6 - 4.58 \times 10^2 \, T^5  + 5.21 \times 10^4 \, T^4  - 3.04 \times 10^6 \, T^3  + 9.61 \times 10^7 \, T^2 - 1.55 \times 10^9 \, T + 9.89 \times 10^9 )$. The RMSE is very small: $0.006$ min. 
To test the validity of the theoretical 6th-order complexity on
$T$, the data were also fitted by polynomials of orders 7, 5, 4, and 3.
Using a 7th-order polynomial does not improve the RMSE.
A 5th-order polynomial fit gives still a very small RMSE of $0.01$ min. (under 1 sec.), while a 4th-order polynomial fit has a small RMSE of $0.03$ min. (about 2 sec.), and a 
3rd-order fit has a much larger RMSE of $0.12$ min.
These results suggest that the predicted runtime of the $\mathrm{RAG}(1, 1)$ algorithm is best fitted by a polynomial of lower order than the 
theoretical 6th-order complexity, with the 4th-order fit $\widehat{y}(T) = 10^7 (6.8 \, T^4  - 1.09 \times 10^3 \, T^2 + 6.7 \times 10^4 \, T - 1.8 \times 10^6)$ appearing to be best.
As for the calibration method with $5$ significant digits, the 2nd- and 3rd-order fits
have roughly equal RMSEs of about $0.02$ min., whereas the 1st-order fit has a large RMSE of $0.22$ min.
Hence, the best fit for the predicted runtime of the calibration method 
is  $O(T^2)$. Yet, note that the $\mathrm{RAG}(1, 1)$
algorithm is faster than the calibration method with 5 significant
digits up to and including horizon $T = 70$.

\begin{figure}[ht!]
\centering
\includegraphics[height=3in,width=6in,keepaspectratio]{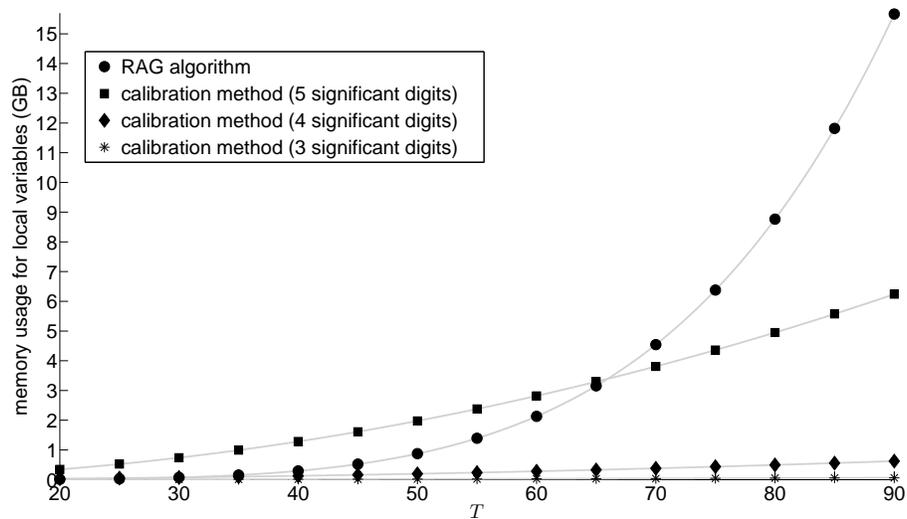}
\caption{Memory usage  vs.\ horizon $T$.}
\label{plot3BernExp}
\end{figure}

As for the actual memory usage of each algorithm
 (for local variables, excluding input and output), 
Figure \ref{plot3BernExp} 
plots 
the required memory storage for floating-point local variables (excluding input and
  output, and using 8-byte double-precision numbers)  versus $T$. 
In the author's implementations, the $\mathrm{RAG}(1, 1)$ algorithm uses $T^5/3 + 2 T^4 + (13/3) T^3 + (15/2) T^2 + (65/6) T + 6$ floating-point storage locations for local variables, whereas the 
calibration method with a grid of size $L$ uses $L (T+1)(T+2)  + L$ storage locations.
Note that, despite its higher complexity, the RAG algorithm
uses less memory than the calibration method with 5 significant digits up to and including 
horizon $T
= 65$.

\begin{figure}[ht!]
\centering
\includegraphics[height=3in,width=6in,keepaspectratio]{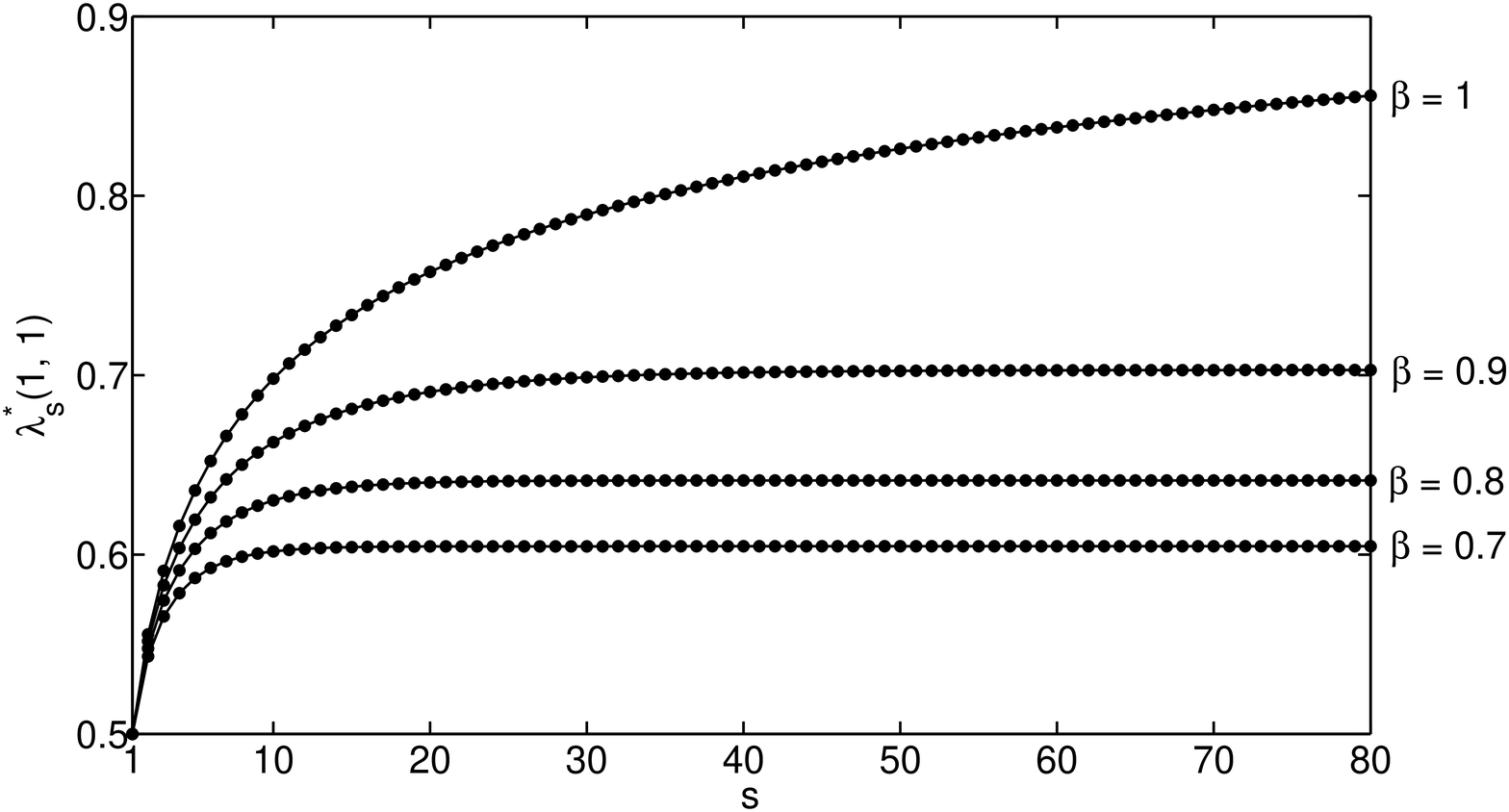}
\caption{Values of index $\lambda^*\big(s, (1, 1)\big)$ vs.\ horizon $s$.}
\label{plot2BernExp}
\end{figure}

Figure \ref{plot2BernExp} plots index $\lambda^*\big(s, (1, 1)\big)$  versus the horizon 
$s = 1, \ldots, 80$, for discount factors $\beta = 0.7, \ldots, 1$.
For each discount factor, the plotted index values were obtained from a single run of algorithm RAG$(1, 1)$ with a horizon of $T = 80$. 

\section{Conclusions}
\label{s:concl}
This paper has introduced a  recursive adaptive-greedy (RAG)
algorithm for the efficient exact computation of a classic index for finite-horizon 
bandits, which performs only arithmetic operations.
The  algorithm has been compared with the standard 
calibration method, which computes approximate index values.
When the latter method is used with a grid having 
the same size as the number of index values to
be evaluated, both methods have the same time and memory complexities. 
Yet, for a fixed grid, the calibration method's time and memory complexity are one order of magnitude lower than those of the RAG algorithm. 
Complementing such theoretical results, 
the computational study reported above shows that,
 if 3 or 4 significant digits of accuracy
suffice, or if the number of states or the horizon are rather large, the calibration method is the best choice. 
Yet, the results also show that
the RAG algorithm outperforms  (with respect to both runtimes and memory)
 the calibration
method with $5$ significant digits  of accuracy in instances of 
moderately large size.

\section*{Acknowledgments}
The author thanks the Associate Editor and an anonymous reviewer
for constructive comments that helped improve the paper, including
suggestions to compare the proposed algorithm with the calibration
method, and to test the validity of the least-squares fits used in the
computational experiments.
This work was 
supported in part by the Spanish Ministry
of Education and Science  under project
 MTM2007-63140 and an I3 faculty endowment grant. 

\bibliographystyle{ijocv081}


\end{document}